\documentclass[twoside]{amsart}

\usepackage{latexsym}
\usepackage{amsthm,amsfonts,amssymb}
\usepackage[leqno]{amsmath}
\usepackage{xspace}
\usepackage[all]{xypic}
\usepackage{hyperref}
\usepackage{cite}

\newcommand{\numberseries}{\mdseries}   %Fontseries used for numbering theorem

\newlength{\thmtopspace}                %Space above theorem
\newlength{\thmbotspace}                %Space below theorem
\newlength{\thmheadspace}               %Space between theorem caption and text
\newlength{\thmindent}                  %For indenting

\newtheoremstyle{bfupright head,slanted body}
                {\thmtopspace}{\thmbotspace}
                {\slshape}{\thmindent}{\bfseries}{.}{\thmheadspace}
                {{\numberseries \thmnumber{(#2) }}\thmnote{#3}}

\newtheoremstyle{bfupright head,upright body}
                {\thmtopspace}{\thmbotspace}
                {\upshape}{\thmindent}{\bfseries}{.}{\thmheadspace}
                {{\numberseries \thmnumber{(#2) }}\thmnote{#3}}

\newtheoremstyle{bfit head,upright body}
                {\thmtopspace}{\thmbotspace}
                {\upshape}{\thmindent}{\upshape}{.}{\thmheadspace}
                {{\numberseries\thmnumber{(#2) }}
                {\bfseries\itshape\thmnote{\negthickspace#3}}}

\newtheoremstyle{it head,upright body}
                {\thmtopspace}{\thmbotspace}
                {\upshape}{\thmindent}{\upshape}{.}{\thmheadspace}
                {{\numberseries\thmnumber{(#2) }}
                {\itshape\thmnote{\negthickspace#3}}}

\newtheoremstyle{fixed bf head,slanted body}
                {\thmtopspace}{\thmbotspace}{\slshape}
                {\thmindent}{\bfseries}{.}{\thmheadspace}
                {{\numberseries \thmnumber{(#2) }}\thmname{#1}\thmnote{ (#3)}}

\newtheoremstyle{fixed bf head,upright body}
                {\thmtopspace}{\thmbotspace}{\upshape}
                {\thmindent}{\bfseries}{.}{\thmheadspace}
                {{\numberseries \thmnumber{(#2) }}\thmname{#1}\thmnote{ (#3)}}

\newtheoremstyle{indented paragraph}
                {\thmtopspace}{\thmbotspace}
                {\upshape}{\thmindent}{\upshape}{}{0pt}
                {\thmnote{#3 }}

\theoremstyle{bfupright head,slanted body}
\newtheorem{res}{}[section]             \newtheorem*{res*}{}

\theoremstyle{bfit head,upright body}
                 \newtheorem*{com*}{}

\theoremstyle{bfupright head,upright body}
\newtheorem{bfhpg}[res]{}               \newtheorem*{bfhpg*}{}

\theoremstyle{it head,upright body}
               \newtheorem*{ithpg*}{}

\theoremstyle{fixed bf head,slanted body}
\newtheorem{thm}[res]{Theorem}          \newtheorem*{thm*}{Theorem}
\newtheorem{prp}[res]{Proposition}      \newtheorem*{prp*}{Proposition}
\newtheorem{cor}[res]{Corollary}        \newtheorem*{cor*}{Corollary}
\newtheorem{lem}[res]{Lemma}            \newtheorem*{lem*}{Lemma}

\theoremstyle{fixed bf head,upright body}
\newtheorem{dfn}[res]{Definition}       \newtheorem*{dfn*}{Definition}
\newtheorem{obs}[res]{Observation}      \newtheorem*{obs*}{Observation}
\newtheorem{rmk}[res]{Remark}           \newtheorem*{rmk*}{Remark}
          \newtheorem*{exa*}{Example}
         \newtheorem*{exe*}{Exercise}
            \newtheorem{stp*}{Setup}
     \newtheorem*{dfns*}{Definitions}
    \newtheorem*{obss*}{Observations}
         \newtheorem*{rmks*}{Remarks}
        \newtheorem*{exas*}{Examples}

\theoremstyle{indented paragraph}
\newtheorem{ipg}{}

\newlength{\thmlistleft}        %leftmargin
\newlength{\thmlistright}       %rightmargin
\newlength{\thmlistpartopsep}   %partopsep
\newlength{\thmlisttopsep}      %topsep
\newlength{\thmlistparsep}      %parsep
\newlength{\thmlistitemsep}     %itemsep

\newcounter{eqc} 
\newenvironment{eqc}{\begin{list}{\upshape (\textit{\roman{eqc}})}%
                    {\usecounter{eqc}%
                        \setlength{\leftmargin}{\thmlistleft}%
                        \setlength{\labelwidth}{\thmlistleft}%
                        \setlength{\rightmargin}{\thmlistright}%
                        \setlength{\partopsep}{\thmlistpartopsep}%
                        \setlength{\topsep}{\thmlisttopsep}%
                        \setlength{\parsep}{\thmlistparsep}%
                        \setlength{\itemsep}{\thmlistitemsep}}}%
                    {\end{list}}%

\newcommand{\eqclbl}[1]{{\upshape(\textit{#1})}}

\newcounter{prt}
\newenvironment{prt}{\begin{list}{\upshape (\alph{prt})}%
                    {\usecounter{prt}%
                        \setlength{\leftmargin}{\thmlistleft}%
                        \setlength{\labelwidth}{\thmlistleft}%
                        \setlength{\rightmargin}{\thmlistright}%
                        \setlength{\partopsep}{\thmlistpartopsep}%
                        \setlength{\topsep}{\thmlisttopsep}%
                        \setlength{\parsep}{\thmlistparsep}%
                        \setlength{\itemsep}{\thmlistitemsep}}}%
                    {\end{list}}%

\newcommand{\prtlbl}[1]{{\upshape(#1)}}

\newcounter{rqm}
\newenvironment{rqm}{\begin{list}{\upshape (\arabic{rqm})}%
                    {\usecounter{rqm}%
                        \setlength{\leftmargin}{\thmlistleft}%
                        \setlength{\labelwidth}{\thmlistleft}%
                        \setlength{\rightmargin}{\thmlistright}%
                        \setlength{\partopsep}{\thmlistpartopsep}%
                        \setlength{\topsep}{\thmlisttopsep}%
                        \setlength{\parsep}{\thmlistparsep}%
                        \setlength{\itemsep}{\thmlistitemsep}}}%
                    {\end{list}}%

                        {\end{list}}%

\newenvironment{itemlist}{\nopagebreak \begin{list}{$\bullet$}%
                       {\setlength{\leftmargin}{\thmlistleft}%
                        \setlength{\labelwidth}{\thmlistleft}%
                        \setlength{\rightmargin}{\thmlistright}%
                        \setlength{\partopsep}{\thmlistpartopsep}%
                        \setlength{\topsep}{\thmlisttopsep}%
                        \setlength{\parsep}{\thmlistparsep}%
                        \setlength{\itemsep}{\thmlistitemsep}}}%
                        {\end{list}}%

                        {\end{list}}%
\newlength{\myindent}
{\setlength{\myindent}{\parindent}\begin{list}{}%
                        {\setlength{\leftmargin}{#1}\setlength{\rightmargin}{#1}%
                        \setlength{\partopsep}{0pt}%
                        \setlength{\topsep}{\thmtopspace}%
                        \setlength{\parsep}{0pt}%
                        \setlength{\itemsep}{0pt}}
                        \item[]}
                        {\end{list}}%

\newenvironment{proof*}{\begin{proof}}{\renewcommand{\qed}{} \end{proof}}

  \newcommand{\step}[1]{$\mathbf{#1^\circ}$}
  \newcommand{\stepref}[1]{$#1^\circ$}

  \newcommand{\proofoftag}[2][:]{(#2)#1}

  \newcommand{\proofofimp}[3][:]{\eqclbl{#2}$\implies$\eqclbl{#3}#1}

\newcommand{\dispand}[1][and]{\hbox to \hsize{#1 \hfill} \nonumber \\}

\newlength{\seqsplit}
\renewcommand{\theequation}{\arabic{equation}}
\numberwithin{equation}{res}

\newcommand{\catbl}{\sqsubset}
\newcommand{\catbr}{\sqsupset}
\newcommand{\catb}{\sqsubset\mspace{-13mu}\sqsupset}

\newcommand{\Cat}[2]{{\sf{#2}}(#1)}

\newcommand{\Catsub}[3]{{\sf{#2}}_{#3}(#1)}
\newcommand{\Catsupsub}[4]{{\sf{#2}}^{\text{\upshape #3}}_{#4}(#1)}

\newcommand{\Catbl}[2]{\Catsub{#1}{#2}{\catbl}}
\newcommand{\Catbr}[2]{\Catsub{#1}{#2}{\catbr}}
\newcommand{\Catb}[2]{\Catsub{#1}{#2}{\catb}}

\newcommand{\Catbrs}[3]{\Catsupsub{#1}{#2}{#3}{\catbr}}
\newcommand{\Catbs}[3]{\Catsupsub{#1}{#2}{#3}{\catb}}

\newcommand{\C}[1][R]{\Cat{#1}{C}}
\newcommand{\Cbl}[1][R]{\Catbl{#1}{C}}
\newcommand{\Cbr}[1][R]{\Catbr{#1}{C}}
\newcommand{\Cb}[1][R]{\Catb{#1}{C}}

\newcommand{\G}[1][R]{\operatorname{G}(#1)}

\newcommand{\A}[1][R]{\Cat{#1}{A}}
\newcommand{\B}[1][R]{\Cat{#1}{B}}

\newcommand{\dd}[2]{{\partial}_{#1}^{#2}}

\renewcommand{\H}[2][\no]{\operatorname{H}_{#1}(#2)}
\newcommand{\amp}[1]{\operatorname{amp}#1}

\newcommand{\cZ}[2]{\operatorname{Z}_{#1}^{#2}}

\newcommand{\cC}[2]{\operatorname{C}_{#1}^{#2}}

\newcommand{\Thr}[2]{#2\mspace{-2mu}\sideset{_{#1}}{}{\operatorname{\sqsupset}}}
\newcommand{\Thl}[2]{\sideset{}{_{#1}}{\operatorname{\sqsubset}}\mspace{-2mu}#2}
\newcommand{\Tsr}[2]{#2\mspace{-2mu}\sideset{_{#1}}{}{\operatorname{\supset}}}
\newcommand{\Tsl}[2]{\sideset{}{_{#1}}{\operatorname{\subset}}\mspace{-2mu}#2}

\newcommand{\Hom}[3][R]{\operatorname{Hom}_{#1}(#2,#3)}
\newcommand{\Ext}[4][R]{\operatorname{Ext}_{#1}^{#2}(#3,#4)}
\newcommand{\tp}[3][R]{#2\otimes_{#1}#3}
\newcommand{\Tor}[4][R]{\operatorname{Tor}^{#1}_{#2}(#3,#4)}

\newcommand{\E}[2][R]{\operatorname{E}_{#1}(#2)}

\newcommand{\D}[1][R]{\Cat{#1}{D}}
\newcommand{\Dbl}[1][R]{\Catbl{#1}{D}}
\newcommand{\Dbr}[1][R]{\Catbr{#1}{D}}
\newcommand{\Db}[1][R]{\Catb{#1}{D}}

\newcommand{\Dfbr}[1][R]{\Catbrs{#1}{D}{f}}
\newcommand{\Dfb}[1][R]{\Catbs{#1}{D}{f}}

\newcommand{\DHom}[3][R]{\operatorname{\mathbf{R}Hom}_{#1}(#2,#3)}
\newcommand{\Dtp}[3][R]{\mbox{\ensuremath{#2\otimes_{#1}^{\mathbf{L}}#3}}}

\newcommand{\RG}[2][\mathfrak{a}]{\mbox{\ensuremath{\mathbf{R}\Gamma_{#1} #2}}}
\newcommand{\LL}[2][\mathfrak{a}]{\mbox{\ensuremath{\mathbf{L}
                                                \Lambda^{#1} #2}}}

\newcommand{\no}{\mspace{-1mu}} %For correct spacing when no ring is specified

\newcommand{\mo}{mor\-phism\xspace}             
\newcommand{\ho}{ho\-mo\-\mo}                   
\newcommand{\qiso}{quasi-isomorphism\xspace}

\newcommand{\eq}{\simeq}
\newcommand{\is}{\cong}

\newcommand{\x}{\pmb{x}}
\newcommand{\defx}[1][t]{\x=x_1,\dots,x_{#1}}

\renewcommand{\a}{\alpha}
\renewcommand{\b}{\beta}

\newcommand{\g}{\gamma}
\renewcommand{\l}{\ell}

\newcommand{\pr}{\pi}

\newcommand{\fe}{\varphi}

\newcommand{\ZZ}{\mathbb{Z}}
\newcommand{\QQ}{\mathbb{Q}}

\newcommand{\m}{\mathfrak{m}}
\newcommand{\n}{\mathfrak{n}}
\newcommand{\p}{\mathfrak{p}}

\newcommand{\Sopp}{S^{\mathrm{opp}}}
\newcommand{\Ropp}{R^{\mathrm{opp}}}
\newcommand{\Rm}{(R,\m)}
\newcommand{\Sn}{(S,\n)}
\newcommand{\Rmk}{(R,\m,k)}

\newcommand{\Ra}{\implies}
\newcommand{\Lra}{\iff}
\newcommand{\lora}{\longrightarrow}
\newcommand{\xla}{\xleftarrow}
\newcommand{\xra}{\xrightarrow}
\newcommand{\xle}{\xla{\;\eq\;}}
\newcommand{\xre}{\xra{\;\eq\;}}
\newcommand{\qeq}{\;=\;}
\newcommand{\mapdef}[4][\longrightarrow]{\mbox{\ensuremath{#2\!: #3 #1 #4}}}

\newcommand{\supremum}[2]{\sup{\{#1\:|\:#2\}}}

\newcommand{\Ker}[1]{\mbox{\ensuremath{\operatorname{Ker}#1}}}
\newcommand{\Coker}[1]{\mbox{\ensuremath{\operatorname{Coker}#1}}}

\renewcommand{\Mc}{\operatorname{\sf Cone}}

\newcommand{\SpecR}{\operatorname{Spec}R}

\newcommand{\dptR}{\operatorname{depth}R}
\newcommand{\dptRp}{\operatorname{depth}R_{\p}}
\newcommand{\dimR}{\operatorname{dim}R}

\newcommand{\fd}[2][R]{\operatorname{fd}_{#1}#2}
\newcommand{\id}[2][R]{\operatorname{id}_{#1}#2}
\newcommand{\pd}[2][R]{\operatorname{pd}_{#1}#2}

\newcommand{\Gdim}[2][R]{\operatorname{G--dim}_{#1}#2}
\newcommand{\Gfd}[2][R]{\operatorname{Gfd}_{#1}#2}
\newcommand{\Gid}[2][R]{\operatorname{Gid}_{#1}#2}
\newcommand{\Gpd}[2][R]{\operatorname{Gpd}_{#1}#2}

\newcommand{\FFD}[1][R]{\operatorname{FFD}(#1)}
\newcommand{\FID}[1][R]{\operatorname{FID}(#1)}
\newcommand{\FPD}[1][R]{\operatorname{FPD}(#1)}

\newcommand{\FGFD}[1][R]{\operatorname{FGFD}(#1)}

\newcommand{\FGPD}[1][R]{\operatorname{FGPD}(#1)}

\newcommand{\wdt}[2][R]{\operatorname{width}_{#1}#2}
\newcommand{\dpt}[2][R]{\operatorname{depth}_{#1}#2}

\newcommand{\ampP}[1]{\amp{(#1)}}
\newcommand{\supP}[1]{\sup{(#1)}}
\newcommand{\infP}[1]{\inf{(#1)}}

\newcommand{\tpP}[3][R]{(\tp[#1]{#2}{#3})}

\newcommand{\DtpP}[3][R]{(\Dtp[#1]{#2}{#3})}

\newcommand{\one}{\ensuremath{\mathord \ast}}
\newcommand{\two}{\ensuremath{\mathord \sharp}}
\newcommand{\three}{\ensuremath{\mathord \star}}
\newcommand{\four}{\ensuremath{\mathord \flat}}

\newcommand{\lhty}[1][D]{\grave{\chi}_#1^{\langle S,R \rangle}}
\newcommand{\rhty}[1][D]{\acute{\chi}_#1^{\langle S,R \rangle}}
\newcommand{\pair}[1][S,R]{\langle #1 \rangle}
\newcommand{\bi}[1]{{}_S#1_R}
\newcommand{\biP}[1]{(\bi{#1})}
\newcommand{\Ca}[1][a]{\operatorname{C}(\mathfrak{#1})}

\renewcommand{\le}{\leqslant}
\renewcommand{\leq}{\leqslant}
\renewcommand{\ge}{\geqslant}
\renewcommand{\geq}{\geqslant}

\begin{document}

\title[ON GORENSTEIN PROJECTIVE, INJECTIVE AND FLAT DIMENSIONS]{On
  Gorenstein projective, injective and \\ flat Dimensions --- A
  functorial description \\ with applications}

\author{Lars Winther Christensen, Anders Frankild, and Henrik Holm}

\address{Lars Winther Christensen, Cryptomathic A/S, Christians Brygge
  28, DK-1559 Copenhagen V, Denmark}

\curraddr{Department of Mathematics, University of Nebraska, Lincoln,
  NE 68588, U.S.A.}

\email{winther@math.unl.edu}

\address{Anders Frankild, Department of Mathematics,
  Universitetsparken 5, DK-2100 Copenhagen \O, Denmark}

\email{frankild@math.ku.dk}

\address{Henrik Holm, Department of Mathematics, Universitetsparken 5,
  DK-2100 Copenhagen \O, Denmark}

\curraddr{Henrik Holm, Department of Mathematical Sciences, University
  of Aarhus, Ny Munkegade Bldg.\ 530, DK-8000 Aarhus C, Denmark}

\email{holm@imf.au.dk}

\date{18 November 2005}



\thanks{L.W.C.\ was partly supported by a grant from the Danish
  Natural Science Research Council}

\thanks{A.F.\ was supported by Lundbeck Fonden, the Danish Natural
  Science Research Council, and Mathematical Sciences Research
  Institute (MSRI)}

\thanks{Part of this work was done at MSRI during the spring semester
  of 2003, when the authors participated in the Program in Commutative
  Algebra.  We thank the institution and program organizers for a very
  stimulating research environment}

\dedicatory{Dedicated to Professor Christian U.~Jensen}

\begin{abstract}
  Gorenstein homological dimensions are refinements of the classical
  homological dimensions, and finiteness singles out modules with
  amenable properties reflecting those of modules over Gorenstein
  rings.
  
  As opposed to their classical counterparts, these dimensions do not
  immediately come with practical and robust criteria for finiteness,
  not even over commutative noetherian local rings. In this paper we
  enlarge the class of rings known to admit good criteria for
  finiteness of Gorenstein dimensions: It now includes, for instance,
  the rings encountered in commutative algebraic geometry and, in the
  non-commutative realm, $k$--algebras with a dualizing complex.
\end{abstract}

\maketitle

%%% SECTION 0
\section*{Introduction}
\label{sec:Introduction}

An important motivation for studying homological dimensions goes back
to 1956 when Auslander, Buchsbaum and Serre proved the following
theorem: A commutative noetherian local ring $R$ is regular if the
residue field $k$ has finite projective dimension and only if all
$R$--modules have finite projective dimension. This introduced the
theme that finiteness of a homological dimension for all modules
singles out rings with special properties.

Subsequent work showed that over any commutative noetherian ring,
modules of finite projective or injective dimension have special
properties resembling those of modules over regular rings. This is one
reason for studying homological dimensions of individual modules.

This paper is concerned with homological dimensions for modules over
associative rings. In the introduction we restrict to a commutative
noetherian local ring $R$ with residue field $k$.
\begin{center}
  $*\ *\ *$
\end{center}
Pursuing the themes described above, Auslander and Bridger
\cite{MAu67,MAuMBr69} introduced a homological dimension designed to
single out modules with properties similar to those of modules over
Gorenstein rings. They called it the G--dimension, and it has strong
parallels to the projective dimension: $R$ is Gorenstein if the
residue field $k$ has finite G--dimension and only if all finitely
generated $R$--modules have finite G--dimension. A finitely generated
$R$--module $M$ of finite G--dimension satisfies an analogue of the
Auslander--Buchsbaum formula:
\begin{equation*}
  \Gdim{M} = \dptR -\dpt{M}.
\end{equation*}
Like other homological dimensions, the G--dimension is introduced by
first defining the modules of dimension $0$, and then using these to
resolve arbitrary modules. Let us recall the definition: A finitely
generated $R$--module $M$ has G--dimension $0$ if
\begin{equation*}
  \Ext{m}{M}{R} = 0 = \Ext{m}{\Hom{M}{R}}{R}\quad \text{for}\:\: m>0
\end{equation*}
and $M$ is reflexive, that is, the canonical map
\begin{equation*}
  M \xra{\quad} \Hom{\Hom{M}{R}}{R}
\end{equation*}
is an isomorphism. Here we encounter a first major difference between
G--dimension and projective dimension. A projective $R$--module $M$ is
described by vanishing of the cohomology functor $\Ext{1}{M}{-}$, and
projectivity of a finitely generated module can even be verified by
computing a single cohomology module, $\Ext{1}{M}{k}$. However,
verification of G--dimension $0$ requires, a priori, the computation
of infinitely many cohomology modules. Indeed, recent work of
Jorgensen and \c{S}ega \cite{DAJLMS} shows that for a reflexive module
$M$ the vanishing of $\Ext{>0}{M}{R}$ and $\Ext{>0}{\Hom{M}{R}}{R}$
cannot be inferred from vanishing of any finite number of these
cohomology modules.

Since the modules of G--dimension $0$ are not described by vanishing
of a (co)homology functor, the standard computational techniques of
homological algebra, like dimension shift, do not effectively apply to
deal with modules of finite G--dimension. This has always been the
Achilles' heel of the theory and explains the interest in finding
alternative criteria for finiteness of this dimension.

G--dimension also differs from projective dimension in that it is
defined only for finitely generated modules. To circumvent this
shortcoming, Enochs and Jenda~\cite{EEnOJn95} proposed to study a
homological dimension based on a larger class of modules: An
$R$--module $M$ is called Gorenstein projective, if there exists an
exact complex
\begin{equation*}
  P = \dots \to P_{1} \xra{\dd{1}{P}}  P_{0} \xra{\dd{0}{P}}  P_{-1}
  \to \cdots
\end{equation*}
of projective modules, such that $M \is \Coker{\dd{1}{P}}$ and
$\Hom{P}{Q}$ is exact for any projective $R$--module $Q$. This
definition allows for non-finitely generated modules and is, indeed,
satisfied by all projective modules. It was already known from
\cite{MAu67} that a finitely generated $R$--module $M$ is of
G--dimension $0$ precisely when there exists an exact complex $L =
\dots \to L_{1} \to L_{0} \to \cdots$ of finitely generated free
modules, such that $M \is \Coker{\dd{1}{L}}$ and $\Hom{L}{R}$ is
exact. Avramov, Buchweitz, Martsinkovsky, and Reiten (see \cite{LWC3})
proved that for finitely generated modules, the Gorenstein projective
dimension agrees with the G--dimension.

Gorenstein flat and injective modules were introduced along the same
lines \cite{EEnOJn95,EJT93}.  Just as the G--dimension has strong
similarities with the projective dimension, the Gorenstein injective
and flat dimensions are, in many respects, similar to the classical
flat and injective dimensions. However, these new dimensions share the
problem encountered already for G--dimension: It is seldom practical
to study modules of finite dimension via modules of dimension 0.

The goal of this paper is to remedy this situation. We do so by
establishing a conjectured characterization of modules of finite
Gorenstein dimensions in terms of vanishing of homology and
invertibility of certain canonical maps. It extends an idea of Foxby,
who found an alternative criterion for finite G--dimension of finitely
generated modules; see \cite{SYs95}.  Before describing the criteria
for finiteness of Gorenstein projective, injective and flat
dimensions, we present a few applications.
\begin{center}
  $*\ *\ *$
\end{center}
The study of Gorenstein dimensions takes cues from the classical
situation; an example: It is a trivial fact that projective modules
are flat but a deep result, due to Gruson--Raynaud \cite{MRnLGr71} and
Jensen \cite{CUJ70}, that flat $R$--modules have finite projective
dimension. For Gorenstein dimensions the situation is more
complicated. It is true but not trivial that Gorenstein projective
$R$--modules are Gorenstein flat; in fact, the proof relies on the
very result by Gruson, Raynaud and Jensen mentioned above. However,
very little is known about Gorenstein projective dimension of
Gorenstein flat $R$--modules. As a first example of what can be gained
from our characterization of modules of finite Gorenstein dimensions,
we address this question; see~\eqref{cor:GorGRJ}:
\begin{res*}[Theorem I]
  If $R$ has a dualizing complex, then the following are equivalent
  for an $R$--module $M$:
  \begin{eqc}
  \item $M$ has finite Gorenstein projective dimension, $\Gpd{M} <
    \infty$.
  \item $M$ has finite Gorenstein flat dimension, $\Gfd{M} < \infty$.
  \end{eqc}
\end{res*}
As the hypothesis of Theorem I indicates, our characterization of
finite Gorenstein dimensions requires the underlying ring to have a
dualizing complex. By Kawasaki's proof of Sharp's conjecture
\cite{TKw02}, this is equivalent to assuming that $R$ is a homomorphic
image of a Gorenstein local ring.

While the Gorenstein analogues of the Auslander--Buchsbaum--Serre
theorem and the Auslander--Buchsbaum formula were among the original
motives for studying G--dimension, the Gorenstein equivalent of
another classic, the Bass formula, has proved more elusive. It was
first established over Gorenstein rings \cite{EEnOJn95a}, later over
Cohen--Macaulay local rings with dualizing module \cite{LWC3}. The
tools invented in this paper enable us to remove the Cohen--Macaulay
hypothesis; see~\eqref{thm:bass}:
\begin{res*}[Theorem II]
  If $R$ has a dualizing complex, and $N$ is a non-zero finitely
  generated $R$--module of finite Gorenstein injective dimension, then
  \begin{equation*}
    \Gid{N} = \dptR.
  \end{equation*}
\end{res*}
As a third application we record the following result, proved in
\eqref{thm:prod_Gflat} and \eqref{thm:Gid_and_limits}:
\begin{res*}[Theorem III]
  If $R$ has a dualizing complex, then any direct product of
  Gorenstein flat $R$--modules is Gorenstein flat, and any direct sum
  of Gorenstein injective modules is Gorenstein injective.
\end{res*}
Over any noetherian ring, a product of flat modules is flat and a sum
of injectives is injective; this is straightforward. The situation for
Gorenstein dimensions is, again, more complicated and, hitherto,
Theorem III was only known for some special rings.
\begin{center}
  $*\ *\ *$
\end{center}
The proofs of Theorems I--III above rely crucially on a description of
finite Gorenstein homological dimensions in terms of two full
subcategories of the derived category of $R$--modules. They are the
so-called Auslander categories, $\A$ and $\B$, associated to the
dualizing complex; they were first studied in \cite{EJX96,HBF94}.

We prove that the modules in $\A$ are precisely those of finite
Gorenstein projective dimension, see theorem~\eqref{thm:mainthm_A},
and the modules in $\B$ are those of finite Gorenstein injective
dimension, see theorem~\eqref{thm:mainthm_B}. For many applications it
is important that these two categories are related through an
equivalence that resembles Morita theory:
\begin{displaymath}
  \xymatrix{\A \ar@<0.7ex>[rrr]^-{D\otimes_R^{\mathbf{L}}-} & {} & {} &
    \B; \ar@<0.7ex>[lll]^-{\DHom{D}{-}}}
\end{displaymath}
where $D$ is the dualizing complex. This may be viewed as an extension
of well-known facts: If $R$ is Cohen--Macaulay, the equivalences above
restrict to the subcategories of modules of finite flat and finite
injective dimension. If $R$ is Gorenstein, that is, $D=R$, the
subcategories of modules of finite flat and finite injective dimension
even coincide, and so do $\A$ and $\B$.

For a Cohen--Macaulay local ring with a dualizing module, this
description of finite Gorenstein homological dimensions in terms of
$\A$ and $\B$ was established in \cite{EJX96}. The present version
extends it in several directions: The underlying ring is not assumed
to be either commutative, or Cohen--Macaulay, or local.

In general, we work over an associative ring with unit. For the main
results, the ring is further assumed to admit a dualizing complex.
Also, we work consistently with complexes of modules. Most proofs, and
even the definition of Auslander categories, require complexes, and it
is natural to state the results in the same generality.
\begin{center}
  $*\ *\ *$
\end{center}
The characterization of finite Gorenstein homological dimensions in
terms of Auslander categories is proved in section~\ref{sec:ABC};
sections 5 and 6 are devoted to applications. The main theorems are
proved through new technical results on preservation of
quasi-isomorphisms; these are treated in section 2.  The first section
fixes notation and prerequisites, and in the third we establish the
basic properties of Gorenstein dimensions in the generality required
for this paper.

%%% SECTION 1
\section{Notation and prerequisites}
\label{sec:notation}

\begin{ipg}
  In this paper, all rings are assumed to be associative with unit,
  and modules are, unless otherwise explicitly stated, left modules.
  For a ring $R$ we denote by $\Ropp$ the opposite ring, and identify
  right $R$--modules with left $\Ropp$--modules in the natural way.
  Only when a module has bistructure, do we include the rings in the
  symbol; e.g., $_S M_R$ means that $M$ is an $(S,\Ropp)$--bimodule.
  
  We consistently use the notation from the appendix of \cite{LWC3}.
  In particular, the category of $R$--complexes is denoted $\C$, and
  we use subscripts $\sqsubset$, $\sqsupset$, and
  $\sqsubset\mspace{-13mu}\sqsupset$ to denote boundedness conditions.
  For example, $\Cbr$ is the full subcategory of $\C$ of right-bounded
  complexes.
  
  The derived category is written $\D$, and we use subscripts
  $\sqsubset$, $\sqsupset$, and $\sqsubset\mspace{-13mu}\sqsupset$ to
  denote homological boundedness conditions.  Superscript ``f''
   signifies that the homology is degreewise finitely generated. Thus,
  $\Dfbr$ denotes the full subcategory of $\D$ of homologically
  right-bounded complexes with finitely generated homology modules.
  The symbol ``$\eq$'' is used to designate isomorphisms in $\D$ and
  quasi-isomorphisms in $\C$.  For the derived category and derived
  functors, the reader is referred to the original texts, Verdier's
  thesis \cite{Verdier} and Hartshorne's notes \cite{RAD}, and further
  to a modern account: the book by Gelfand and Manin \cite{GM}.
\end{ipg}  

Next, we review a few technical notions for later use.

\begin{dfn}
  \label{dfn:dc}
  Let $S$ and $R$ be rings. If $S$ is left noetherian and $R$ is right
  noetherian, we refer to the ordered pair $\pair$ as a
  \emph{noetherian pair} of rings.
  
  A \emph{dualizing complex} for a noetherian pair of rings $\pair$ is
  a complex $\bi{D}$ of bimodules meeting the requirements:
  \begin{rqm}
  \item The homology of $D$ is bounded and degreewise finitely
    generated over $S$ and over $\Ropp$.
  \item There exists a quasi-isomorphism of complexes of bimodules,
    $\bi{P} \xre \bi{D}$, where $\bi{P}$ is right-bounded and consists
    of modules projective over both $S$ and $\Ropp$.
  \item There exists a quasi-isomorphism of complexes of bimodules,
    $\bi{D} \xre \bi{I}$, where $\bi{I}$ is bounded and consists of
    modules injective over both $S$ and $\Ropp$.
  \item The homothety morphisms
    \begin{align*}
      \lhty \colon _SS_S &\longrightarrow
      \DHom[\Ropp]{\bi{D}}{\bi{D}}\quad \text{and}\\
      \rhty \colon _RR_R &\longrightarrow \DHom[S]{\bi{D}}{\bi{D}},
    \end{align*}
    are bijective in homology. That is to say,
    \begin{itemlist}
    \item $\lhty$ is invertible in $\D[S]$ (equivalently, invertible
      in $\D[\Sopp]$), and
    \item $\rhty$ is invertible in $\D$ (equivalently, invertible in
      $\D[\Ropp]$).
    \end{itemlist}
  \end{rqm}
  If $R$ is both left and right noetherian (e.g., commutative and
  noetherian), then a dualizing complex for $R$ means a dualizing
  complex for the pair $\pair[R,R]$ (in the commutative case the two
  copies of $R$ are tacitly assumed to have the same action on the
  modules).
\end{dfn}

For remarks on this definition and comparison to other notions of
dualizing complexes in non-commutative algebra, we refer to the
appendix. At this point we just want to mention that \eqref{dfn:dc} is
a natural extension of existing definitions: When $\pair$ is a
noetherian pair of algebras over a field, definition \eqref{dfn:dc}
agrees with the one given by Yekutieli--Zhang
\cite{Yekutieli-Zhang:RADC}. If $R$ is commutative and noetherian,
then \eqref{dfn:dc} is clearly the same as Grothendieck's definition
\cite[V.\S 2]{RAD}.

The next result is proved in the appendix, see page~\pageref{proof:OppPair}.

\begin{prp}
  \label{prp:OppPair}
  Let $\pair$ be a noetherian pair. A complex $\bi{D}$ is dualizing
  for $\pair$ if and only if it is dualizing for the pair
  $\pair[\Ropp,\Sopp]$.
\end{prp}

\begin{bfhpg}[Equivalence]
  \label{bfhpg:Foxby_eq}
  If $\bi{D}$ is a dualizing complex for the noetherian pair $\pair$,
  we consider the adjoint pairs of functors,
  \begin{displaymath}
    \xymatrix{\D \ar@<0.5ex>[rrr]^-{\bi{D}\otimes_R^{\bf L}-} & {} & {} & 
      \D[S]. \ar@<0.5ex>[lll]^-{\DHom[S]{\bi{D}}{-}} 
    } 
  \end{displaymath}
  These functors are represented by $\tp{\bi{P}}{-}$ and
  $\Hom[S]{\bi{P}}{-}$, where $\bi{P}$ is as in \eqref{dfn:dc}(2); see
  also the appendix.
  
  The \emph{Auslander categories} $\A$ and $\B[S]$ with respect to the
  dualizing complex $\bi{D}$ are defined in terms of natural
  transformations being isomorphisms:
  \begin{equation*}
    \A = \left\{ X \in \Db \: 
      \left|
        \begin{array}{l}
          \mbox{ $\eta_X \colon X \xre \DHom[S]{\bi{D}}{\Dtp{\bi{D}}{X}}$ is an iso-} \\
          \mbox{ morphism in $\D$, and $\Dtp{\bi{D}}{X}$ is bounded} 
        \end{array}
      \right.
    \right\},
  \end{equation*}
  and
  \begin{equation*}
    \B[S] = \left\{ Y \in \Db[S] \: 
      \left|
        \begin{array}{l}
          \mbox{ $\varepsilon_Y \colon \Dtp{\bi{D}}{\DHom[S]{\bi{D}}{Y}}\xre Y$ is an
            isomor-} \\ 
          \mbox{ phism in $\D[S]$, and $\DHom[S]{\bi{D}}{Y}$ is bounded}
        \end{array}
      \right.
    \right\}.
  \end{equation*}
  All $R$--complexes of finite flat dimension belong to $\A$, while
  $S$--complexes of finite injective dimension belong to $\B[S]$,
  cf.~\cite[thm.~(3.2)]{LLAHBF97}.
  
  The Auslander categories $\A$ and $\B[S]$ are clearly triangulated
  subcategories of $\D$ and $\D[S]$, respectively, and the adjoint
  pair $(\Dtp{\bi{D}}{-}, \DHom[S]{\bi{D}}{-})$ restricts to an
  equivalence:
  \begin{displaymath}
    \xymatrix{\A
      \ar@<0.7ex>[rrr]^-{\bi{D}\otimes_R^{\mathrm{\mathbf{L}}}-} 
      & {} & {} &
      \B[S]. \ar@<0.7ex>[lll]^-{\DHom[S]{\bi{D}}{-}}
    }
  \end{displaymath}
  In the commutative setting, this equivalence, introduced in
  \cite{LLAHBF97}, is sometimes called \emph{Foxby equivalence}.
\end{bfhpg}

\begin{bfhpg}[Finitistic dimensions]
  \label{finitistic dimensions}
  We write $\FPD$ for the (left) finitistic projective dimension of
  $R$, i.e.
  \begin{displaymath}
    \FPD = \sup \left\{ \, \pd{M} \,
      \left|
        \begin{array}{l}
          \mbox{$M$ is an $R$--module of}\\
          \mbox{finite projective dimension}
        \end{array}
      \right.
    \right\}.
  \end{displaymath}
  Similarly, we write $\FID$ and $\FFD$ for the (left) finitistic
  injective and (left) finitistic flat dimension of $R$.
  
  When $R$ is commutative and noetherian, it is well-known from
  \cite[cor.~5.5]{HBs62} and \cite[II.~thm.~(3.2.6)]{MRnLGr71} that
  \begin{align*}
    \FID \, = \, \FFD \, \le \, \FPD \, = \, \dimR.
  \end{align*}
  If, in addition, $R$ has a dualizing complex, then $\dimR$ is finite
  by \cite[cor.~V.7.2]{RAD}, and hence so are the finitistic
  dimensions.
\end{bfhpg}

Jensen \cite[prop.~6]{CUJ70} proved that if $\FPD$ is finite, then any
$R$--module of finite flat dimension has finite projective dimension
as well. For non-commutative rings, finiteness of $\FPD$ has proved
difficult to establish; indeed, even for finite dimensional algebras
over a field it remains a conjecture. Therefore, the following result
is of interest:

\begin{prp}
  \label{prp:PJresult}
  Assume that the noetherian pair $\pair$ has a dualizing complex
  $\bi{D}$. If $X \in \D$ has finite $\fd{X}$, then there is an
  inequality,
  \begin{displaymath}  
    \pd{X} \le \max\big\{\id[S]{\biP{D}} + \supP{\Dtp{\bi{D}}{X}} \,,\,
    \sup{X}\big\} < \infty. 
  \end{displaymath}
  Moreover, $\FPD$ is finite if and only if $\,\FFD$ is finite.
\end{prp} 

\begin{proof*}
  When $R=S$ is commutative, this was proved by Foxby
  \cite[cor.~3.4]{HBF77}. Recently, J{\o}rgensen \cite{PJ:Preprint}
  has generalized the proof to the situation where $R$ and $S$ are
  $k$--algebras. The further generalization stated above is proved in
  \eqref{bfhpg:PJproof}.
\end{proof*}

Dualizing complexes have excellent duality properties; we shall only
need that for commutative rings:

\begin{bfhpg}[Duality]
  \label{bfhpg:dagger duality}
  Let $R$ be commutative and noetherian with a dualizing complex $D$.
  Grothendieck \cite[V.\S 2]{RAD} considered the functor $-^{\dagger}
  = \DHom{-}{D}$. As noted ibid.\ $-^{\dagger}$ sends $\Dfb$ to itself
  and, in fact, gives a duality on that category.  That is, there is
  an isomorphism
  \begin{equation}
    \label{equ:dagger duality}
    X \xre X^{\dagger \, \dagger} = \DHom{\DHom{X}{D}}{D}
  \end{equation}
  for $X\in\Dfb$, see~\cite[prop.~V.2.1]{RAD}.
\end{bfhpg}

\begin{ipg}
  We close this section by recalling the definitions of Gorenstein
  homological dimensions; they go back to
  \cite{EJT93,EEnOJn93,EEnOJn95,LWC3}.
\end{ipg}

\begin{bfhpg}[Gorenstein projective dimension]
  \label{dfn:gproj}
  An $R$--module $A$ is \emph{Gorenstein projective} if there exists
  an exact complex $P$ of projective modules, such that $A$ is
  isomorphic to a cokernel of $P$, and $\H{\Hom{P}{Q}}=0$ for all
  projective $R$--modules $Q$. Such a complex $P$ is called a
  \emph{complete projective resolution} of $A$.
  
  The \emph{Gorenstein projective dimension}, $\Gpd{X}$, of $X\in\Dbr$
  is defined as
  \begin{align*}
    \Gpd{X} = \inf \left\{ \, \sup \{ \l \in \ZZ \, | \, A_{\l} \ne 0
      \} \, \left|
        \begin{array}{l}
          \mbox{$A\in \Cbr$ is isomorphic to $X$ in $\D$} \\
          \mbox{and every $A_{\l}$ is Gorenstein projective}
        \end{array}
      \right.  \right\}.
  \end{align*}
\end{bfhpg}

\begin{bfhpg}[Gorenstein injective dimension]
  \label{dfn:gim}
  The definitions of \emph{Gorenstein injective modules} and
  \emph{complete injective resolutions} are dual to the ones given in
  \eqref{dfn:gproj}, see also \cite[(6.1.1) and (6.2.2)]{LWC3}. The
  \emph{Gorenstein injective dimension}, $\Gid{Y}$, of $Y\in\Dbl$ is
  defined as
  \begin{displaymath}
    \Gid{Y} = \inf \left\{ \, \sup \{ \l \in \ZZ \, | \, B_{-\l} \ne 0 \} \,
      \left|
        \begin{array}{l}
          \mbox{$B\in \Cbl$ is isomorphic to $Y$ in $\D$} \\
          \mbox{and every $B_{\l}$ is Gorenstein injective}
        \end{array}
      \right.
    \right\}.
  \end{displaymath}
\end{bfhpg}

\begin{bfhpg}[Gorenstein flat dimension]
  \label{dfn:gflat}
  An $R$--module $A$ is \emph{Gorenstein flat} if there exists an
  exact complex $F$ of flat modules, such that $A$ is isomorphic to a
  cokernel of $F$, and $\H{\tp{J}{F}}=0$ for all injective
  $\Ropp$--modules $J$. Such a complex $F$ is called a \emph{complete flat
    resolution} of $A$.
  
  The definition of the \emph{Gorenstein flat dimension}, $\Gfd{X}$,
  of $X\in\Dbr$ is similar to that of the Gorenstein projective
  dimension given in \eqref{dfn:gproj}, see also \cite[(5.2.3)]{LWC3}.
\end{bfhpg}

%%% SECTION 2
\section{Ubiquity of quasi-isomorphisms}
\label{sec:Ubiquity}

\begin{ipg}
  In this section we establish some important, technical results on
  preservation of quasi-isomorphisms. It is, e.g., a crucial
  ingredient in the proof of the main theorem \eqref{thm:mainthm_A}
  that the functor $\tp{-}{A}$ preserves certain quasi-isomorphisms,
  when $A$ is a Gorenstein flat module. This is established in theorem
  \eqref{thm:P->I} below. An immediate consequence of this result is
  that Gorenstein flat modules may sometimes replace real flat modules
  in representations of derived tensor products. This corollary,
  \eqref{prp:rfst-CovDHom}, plays an important part in the proof of
  theorem \eqref{thm:GFDT}.
  
  Similar results on representations of the derived Hom functor are
  used in the proofs of theorems~\eqref{thm:GPDT} and
  \eqref{thm:GIDT}. These are also established below.
  
  The section closes with three approximation results for modules of
  finite Gorenstein homological dimension. We recommend that this
  section is consulted as needed rather than read linearly.
\end{ipg}

The first lemmas are easy consequences of the definitions of
Gorenstein projective, injective, and flat modules.

\begin{lem}
  \label{lem:gp-vanish}
  If $M$ is a Gorenstein projective $R$--module, then
  $\Ext{m}{M}{T}=0$ for all $m>0$ and all $R$--modules $T$ of finite
  projective or finite injective dimension.
\end{lem}

\begin{proof}
  For a module $T$ of finite projective dimension, the vanishing of
  $\Ext{m}{M}{T}$ is an immediate consequence of the definition of
  Gorenstein projective modules.
  
  Assume that $\id{T}= n< \infty$. Since $M$ is Gorenstein projective,
  we have an exact sequence,
  \begin{align*}
    0 \to M \to P_0 \to P_{-1} \to \dots \to P_{1-n} \to C \to 0,
  \end{align*}
  where the $P$'s are projective modules. Breaking this sequence into
  short exact ones, we see that $\Ext{m}{M}{T} = \Ext{m+n}{C}{T}$ for
  $m>0$, so the Exts vanish as desired since $\Ext{w}{-}{T}=0$ for
  $w>n$.
\end{proof}

Similarly one establishes the next two lemmas.

\begin{lem}
  \label{lem:gi-vanish}
  If $N$ is a Gorenstein injective $R$--module, then $\Ext{m}{T}{N}=0$
  for all $m>0$ and all $R$--modules $T$ of finite projective or
  finite injective dimension. \qed
\end{lem}

\begin{lem}
  \label{lem:gf-vanish}
  If $M$ is a Gorenstein flat $R$--module, then $\Tor{m}{T}{M}=0$ for
  all $m>0$ and all $\Ropp$--modules $T$ of finite flat or finite
  injective dimension. \qed
\end{lem}

From lemma~\eqref{lem:gp-vanish} it is now a three step process to
arrive at the desired results on preservation of quasi-isomorphisms by
the Hom functor. We give proofs for the results regarding the
covariant Hom functor; those on the contravariant functor have similar
proofs.

\begin{lem}
  \label{Hom-Lemma}
  Assume that $X, Y \in \C$ with either $X \in \Cbr$ or $Y \in \Cbr$.
  If $\,\H{\Hom{X_\l}{Y}} = 0$ for all $\l \in \ZZ$, then
  $\H{\Hom{X}{Y}} = 0$.
  
  {\rm This result can be found in \cite[lem.~(6.7)]{HHA}. However,
    since this reference is not easily accessible, we provide an
    argument here:}
\end{lem}

\begin{proof}
  Fix an $n \in \ZZ$; we shall prove that $\H[n]{\Hom{X}{Y}}=0$.
  Since
  \begin{displaymath}
    \H[n]{\Hom{X}{Y}}=\H[0]{\Hom{X}{\Sigma^{-n}Y}}
  \end{displaymath}
  we need to show that every morphism $\alpha \colon X \lora
  \Sigma^{-n}Y$ is null-homotopic. Thus for a given morphism $\alpha$
  we must construct a family $(\gamma_m)_{m\in\ZZ}$ of degree 1 maps,
  $\gamma_m \colon X_m \lora (\Sigma^{-n}Y)_{m+1} = Y_{n+m+1}$, such
  that
  \begin{displaymath}
    \tag*{\ensuremath{(\one)_m}}
    \alpha_m = \gamma_{m-1}\partial_m^X +
    \partial_{m+1}^{\Sigma^{-n}Y}\gamma_m
  \end{displaymath}
  for all $m\in\ZZ$. We do so by induction on $m$: Since $X$ or
  $\Sigma^{-n}Y$ is in $\Cbr$ we must have $\gamma_m=0$ for $m \ll 0$.
  For the inductive step, assume that $\gamma_m$ has been constructed
  for all $m<\widetilde{m}$.  By assumption
  $\H[]{\Hom{X_{\widetilde{m}}}{\Sigma^{-n}Y}}=0$, so applying
  $\Hom{X_{\widetilde{m}}}{-}$ to
  \begin{displaymath}
    \Sigma^{-n}Y = \ \cdots \xra{} Y_{n+\widetilde{m}+1}
    \xra{\partial_{\widetilde{m}+1}^{\Sigma^{-n}Y}} Y_{n+\widetilde{m}}
    \xra{\partial_{\widetilde{m}}^{\Sigma^{-n}Y}} Y_{n+\widetilde{m}-1} \xra{}
    \cdots 
  \end{displaymath}
  yields an exact complex. Using that $\alpha$ is a morphism and that
  \ensuremath{(\one)_{\widetilde{m}-1}} holds, we see that
  $\alpha_{\widetilde{m}}-\gamma_{\widetilde{m}-1}\partial_{\widetilde{m}}^X$
  is in the kernel of
  $\Hom{X_{\widetilde{m}}}{\partial_{\widetilde{m}}^{\Sigma^{-n}Y}}$:
  \begin{align*}
    \partial_{\widetilde{m}}^{\Sigma^{-n}Y}
    (\alpha_{\widetilde{m}}-\gamma_{\widetilde{m}-1}\partial_{\widetilde{m}}^X)
    &= (\alpha_{\widetilde{m}-1} -
    \partial_{\widetilde{m}}^{\Sigma^{-n}Y}\gamma_{\widetilde{m}-1})
    \partial_{\widetilde{m}}^X \\
    &=
    (\gamma_{\widetilde{m}-2}\partial_{\widetilde{m}-1}^X)\partial_{\widetilde{m}}^X
    \\ &= 0.
  \end{align*}
  By exactness,
  $\alpha_{\widetilde{m}}-\gamma_{\widetilde{m}-1}\partial_{\widetilde{m}}^X$
  is also in the image of
  $\Hom{X_{\widetilde{m}}}{\partial_{\widetilde{m}+1}^{\Sigma^{-n}Y}}$,
  which means that there exists $\gamma_{\widetilde{m}} \colon
  X_{\widetilde{m}} \lora Y_{n+\widetilde{m}+1}$ such that
  $\partial_{\widetilde{m}+1}^{\Sigma^{-n}Y} \gamma_{\widetilde{m}} =
  \alpha_{\widetilde{m}}-\gamma_{\widetilde{m}-1}\partial_{\widetilde{m}}^X$.
\end{proof}

\begin{lem}
  \label{Dualt Hom-Lemma}
  Assume that $X, Y \in \C$ with either $X \in \Cbl$ or $Y \in \Cbl$.
  If $\,\H{\Hom{X}{Y_\l}} = 0$ for all $\l \in \ZZ$, then
  $\H{\Hom{X}{Y}} = 0$.\qed
\end{lem}

\begin{prp}
  \label{prp:vp_hom_cov}
  Consider a class $\mathfrak{U}$ of $R$--modules, and let $\alpha
  \colon X \lora Y$ be a morphism in $\C$, such that
  \begin{gather*}
    \Hom{U}{\alpha} \colon \Hom{U}{X} \xre \Hom{U}{Y}
  \end{gather*}
  is a quasi-isomorphism for every module $U \in \mathfrak{U}$.  Let
  $\widetilde{U} \in \C$ be a complex consisting of modules from
  $\mathfrak{U}$. The induced morphism,
  \begin{gather*}
    \Hom{\widetilde{U}}{\alpha} \colon \Hom{\widetilde{U}}{X} \lora
    \Hom{\widetilde{U}}{Y},
  \end{gather*}
  is then a quasi-isomorphism, provided that either
  \begin{prt}
  \item $\widetilde{U}\in\Cbr$ or
  \item $X,Y\in\Cbr$.
  \end{prt}
\end{prp}

\begin{proof}
  Under either hypothesis, \prtlbl{a} or \prtlbl{b}, we must verify
  exactness of
  \begin{gather*}
    \Mc({\Hom{\widetilde{U}}{\alpha}}) \, \is \,
    \Hom{\widetilde{U}}{\Mc(\alpha)}.
  \end{gather*}
  Condition \prtlbl{b} implies that $\Mc(\alpha) \in \Cbr$. In any
  event, lemma \eqref{Hom-Lemma} informs us that it suffices to show
  that the complex $\Hom{\widetilde{U}_\l}{\Mc(\alpha)}$ is exact for
  all $\l \in \ZZ$, and this follows as all
  \begin{gather*}
    \Hom{\widetilde{U}_\l}{\alpha} \colon \Hom{\widetilde{U}_\l}{X}
    \xre \Hom{\widetilde{U}_\l}{Y}
  \end{gather*}
  are assumed to be quasi-isomorphisms in $\C$.
\end{proof}

\begin{prp}
  \label{prp:vp_hom_contra}
  Consider a class $\mathfrak{V}$ of $R$--modules, and let $\alpha
  \colon X \lora Y$ be a morphism in $\C$, such that
  \begin{gather*}
    \Hom{\alpha}{V} \colon \Hom{Y}{V} \xre \Hom{X}{V}
  \end{gather*}
  is a quasi-isomorphism for every module $V \in \mathfrak{V}$.  Let
  $\widetilde{V} \in \C$ be a complex consisting of modules from
  $\mathfrak{V}$. The induced morphism,
  \begin{gather*}
    \Hom{\alpha}{\widetilde{V}} \colon \Hom{Y}{\widetilde{V}} \lora
    \Hom{X}{\widetilde{V}},
  \end{gather*}
  is then a quasi-isomorphism, provided that either
  \begin{prt}
  \item $\widetilde{V}\in\Cbl$ or
  \item $X,Y\in\Cbl$. \qed
  \end{prt}
\end{prp}

\begin{thm}
  \label{thm:Hom_quiso1}
  Let $V \xre W$ be a quasi-isomorphism between $R$--complexes, where
  each module in $V$ and $W$ has finite projective dimension or finite
  injective dimension.  If $A\in\Cbr$ is a complex of Gorenstein
  projective modules, then the induced morphism
  \begin{equation*}
    \Hom{A}{V} \lora \Hom{A}{W}
  \end{equation*}
  is a quasi-isomorphism under each of the next two conditions.
  \begin{prt}
  \item $V,W\in\Cbl$ or
  \item $V,W\in\Cbr$.
  \end{prt}
\end{thm}

\begin{proof}
  By proposition \eqref{prp:vp_hom_cov}\prtlbl{a} we may immediately
  reduce to the case, where $A$ is a Gorenstein projective module. In
  this case we have quasi-isomorphisms $\mu \colon P \xre A$ and $\nu
  \colon A \xre \widetilde{P}$ in $\C$, where $P \in \Cbr$ and
  $\widetilde{P}\in \Cbl$ are, respectively, the ``left half'' and
  ``right half'' of a complete projective resolution of $A$.
  
  Let $T$ be any $R$--module of finite projective or finite injective
  dimension. Lemma \eqref{lem:gp-vanish} implies that a complete
  projective resolution stays exact when the functor $\Hom{-}{T}$ is
  applied to it. In particular, the induced morphisms
  \begin{align*}
    \tag{\one} \Hom{\mu}{T} \colon \Hom{A}{T} &\xre \Hom{P}{T},\quad\text{and}\\
    \tag{\two} \Hom{\nu}{T} \colon \Hom{\widetilde{P}}{T} &\xre
    \Hom{A}{T}
  \end{align*}
  are quasi-isomorphisms.  From (\one) and proposition
  \eqref{prp:vp_hom_contra}\prtlbl{a} it follows that under assumption
  \prtlbl{a} both $\Hom{\mu}{V}$ and $\Hom{\mu}{W}$ are
  quasi-isomorphisms. In the commutative diagram
  \begin{gather*}
    \xymatrix{\Hom{A}{V} \ar[d]_-{\Hom[\no]{\mu}{V}}^-{\simeq} \ar[r]
      & \Hom{A}{W} \ar[d]^-{\Hom[\no]{\mu}{W}}_-{\simeq} \\
      \Hom{P}{V} \ar[r]^-{\simeq} & \Hom{P}{W} }
  \end{gather*}
  the lower horizontal morphism is obviously a quasi-isomorphism, and
  this makes the induced morphism $\Hom{A}{V} \longrightarrow
  \Hom{A}{W}$ a quasi-isomorphism as well.
  
  Under assumption \prtlbl{b}, the induced morphism
  $\Hom{\widetilde{P}}{V} \lora \Hom{\widetilde{P}}{W}$ is a
  quasi-isomorphism by proposition \eqref{prp:vp_hom_cov}\prtlbl{b}.
  As the induced morphisms (\two) are quasi-isomorphisms, it follows
  by proposition \eqref{prp:vp_hom_contra}\prtlbl{b} that so are
  $\Hom{\nu}{V}$ and $\Hom{\nu}{W}$. From the commutative diagram
  \begin{gather*}
    \xymatrix{\Hom{A}{V} \ar[r] &
      \Hom{A}{W}  \\
      \Hom{\widetilde{P}}{V} \ar[u]^-{\Hom[\no]{\nu}{V}}_-{\simeq}
      \ar[r]^-{\simeq} & \Hom{\widetilde{P}}{W}
      \ar[u]_-{\Hom[\no]{\nu}{W}}^-{\simeq} }
  \end{gather*}
  we conclude that also its top vertical morphism is a
  quasi-isomorphism.
\end{proof}

\begin{thm}
  \label{thm:Hom_quiso2}
  Let $V \xre W$ be a quasi-isomorphism between $R$--complexes, where
  each module in $V$ and $W$ has finite projective dimension or finite
  injective dimension.  If $B\in\Cbl$ is a complex of Gorenstein
  injective modules, then the induced morphism
  \begin{equation*}
    \Hom{W}{B} \lora \Hom{V}{B}
  \end{equation*}
  is a quasi-isomorphism under each of the next two conditions.
  \begin{prt}
  \item $V,W\in\Cbr$,
  \item $V,W\in\Cbl$ or \qed
  \end{prt}
\end{thm}

\begin{cor}
  \label{prp:rpst-CovDHom}
  Assume that $X\eq A$, where $A\in\Cbr$ is a complex of Gorenstein
  projective modules. If $\,U\eq V$, where $V\in\Cbl$ is a complex in
  which each module has finite projective dimension or finite
  injective dimension, then
  \begin{equation*}
    \DHom{X}{U}\eq \Hom{A}{V}.
  \end{equation*}
\end{cor}

\begin{proof}
  We represent $\DHom{X}{U} \eq \DHom{A}{V}$ by the complex
  $\Hom{A}{I}$, where $V\xre I\in\Cbl$ is an injective resolution of
  $V$.  From \eqref{thm:Hom_quiso1}\prtlbl{a} we get a
  quasi-isomorphism $\Hom{A}{V} \xre \Hom{A}{I}$, and the result
  follows.
\end{proof}

\begin{rmk}
  There is a variant of theorem \eqref{thm:Hom_quiso1} and corollary
  \eqref{prp:rpst-CovDHom}. If $R$ is commutative and noetherian, and
  $A$ is a complex of \textsl{finitely generated} Gorenstein projective
  $R$--modules, then we may relax the requirements on the modules in
  $V$ and $W$ without changing the conclusions of
  \eqref{thm:Hom_quiso1} and \eqref{prp:rpst-CovDHom}: It is
  sufficient that each module in $V$ and $W$ has finite \textsl{flat}
  or finite injective dimension. This follows immediately from the
  proofs of \eqref{thm:Hom_quiso1} and \eqref{prp:rpst-CovDHom}, when
  one takes \cite[prop.~(4.1.3)]{LWC3} into account.
\end{rmk}

\begin{cor}
  \label{prp:rist-CovDHom}
  Assume that $Y\eq B$, where $B\in\Cbl$ is a complex of Gorenstein
  injective modules. If $\,U\eq V$, where $V\in\Cbr$ is a complex in
  which each module has finite projective dimension or finite
  injective dimension, then
  \begin{xxalignat}{3}
    && \DHom{U}{Y} &\eq \Hom{V}{B}.  &&\qed
  \end{xxalignat}
\end{cor}

Next, we turn to tensor products and Gorenstein flat modules.  The
first lemma follows by applying Pontryagin duality to lemma
\eqref{Hom-Lemma} for $\Ropp$.

\begin{lem}
  \label{tensor-Lemma}
  Assume that $X \in \C[\Ropp]$ and $Y \in \C$ with either $X \in
  \Cbr[\Ropp]$ or $Y \in \Cbl$. If $\,\H{\tp{X_\l}{Y}} = 0$ for all $\l
  \in \ZZ$, then $\H {\tp{X}{Y}} = 0$.\qed
\end{lem}

\begin{prp}
  \label{prp:vp}
  Consider a class $\mathfrak{W}$ of $\Ropp$--modules, and let $\alpha
  \colon X \lora Y$ be a morphism in $\C$, such that
  \begin{gather*}
    \tp{W}{\alpha} \colon \tp{W}{X} \xre \tp{W}{Y}
  \end{gather*}
  is a quasi-isomorphism for every module $W \in \mathfrak{W}$.  Let
  $\widetilde{W} \in \C[\Ropp]$ be a complex consisting of modules
  from $\mathfrak{W}$. The induced morphism,
  \begin{gather*}
    \tp{\widetilde{W}}{\alpha} \colon \tp{\widetilde{W}}{X} \lora
    \tp{\widetilde{W}}{Y},
  \end{gather*}
  is then a quasi-isomorphism, provided that either
  \begin{prt}
  \item $\widetilde{W}\in\Cbr[\Ropp]$ or
  \item $X,Y\in\Cbl$.
  \end{prt}
\end{prp}

\begin{proof}
  It follows by lemma~\eqref{tensor-Lemma} that
  $\Mc{(\tp{\widetilde{W}}{\alpha})} \eq
  \tp{\widetilde{W}}{\Mc{(\alpha})}$ is exact under either assumption,
  \prtlbl{a} or \prtlbl{b}.
\end{proof}

\begin{thm}
  \label{thm:P->I}
  Let $V \xre W$ be a quasi-isomorphism between complexes of
  $\Ropp$--modules, where each module in $V$ and $W$ has finite
  injective dimension or finite flat dimension.  If $A \in \Cbr$ is a
  complex of Gorenstein flat modules, then the induced morphism
  \begin{equation*}
    \tp{V}{A} \lora \tp{W}{A}
  \end{equation*}
  is a quasi-isomorphism under each of the next two conditions.
  \begin{prt}
  \item $V,W\in\Cbr[\Ropp]$ or
  \item $V,W\in\Cbl[\Ropp]$.
  \end{prt}
\end{thm}

\begin{proof}
  Using proposition \eqref{prp:vp}\prtlbl{a}, applied to $\Ropp$, we
  immediately reduce to the case, where $A$ is a Gorenstein flat
  module. In this case we have quasi-isomorphisms $\mu \colon F \xre
  A$ and $\nu \colon A \xre \widetilde{F}$ in $\C$, where $F \in \Cbr$
  and $\widetilde{F}\in \Cbl$ are complexes of flat modules. To be
  precise, $F$ and $\widetilde{F}$ are, respectively, the ``left
  half'' and ``right half'' of a complete flat resolution of $A$. The
  proof now continues as the proof of theorem~\eqref{thm:Hom_quiso1};
  only using proposition~\eqref{prp:vp} instead of
  \eqref{prp:vp_hom_cov} and \eqref{prp:vp_hom_contra}, and
  lemma~\eqref{lem:gf-vanish} instead of \eqref{lem:gp-vanish}.
\end{proof}

\begin{cor}
  \label{prp:rfst-CovDHom}
  Assume that $X\eq A$, where $A\in\Cbr$ is a complex of Gorenstein
  flat modules. If $\,U\eq V$, where $V\in \Cbr[\Ropp]$ is a complex in
  which each module has finite flat dimension or finite injective
  dimension, then
  \begin{equation*}
    \Dtp{U}{X}\eq \tp{V}{A}.
  \end{equation*}
\end{cor}

\begin{proof}
  We represent $\Dtp{U}{X} \eq \Dtp{V}{A}$ by the complex $\tp{P}{A}$,
  where\linebreak[4] \mbox{$\Cbr[\Ropp] \ni P \xre V$} is a projective
  resolution of $V$.  By theorem~\eqref{thm:P->I}\prtlbl{a} we get a
  quasi-isomorphism $\tp{P}{A}\xre \tp{V}{A}$, and the desired result
  follows.
\end{proof}

\begin{ipg}
  The Gorenstein dimensions refine the classical homological
  dimensions. On the other hand, the next three lemmas show that a
  module of finite Gorenstein projective/injective/flat dimension can
  be approximated by a module, for which the corresponding classical
  homological dimension is finite.
\end{ipg}

\begin{lem}
  \label{lem:helping_sequence_for_GPD}
  Let $M$ be an $R$--module of finite Gorenstein projective dimension.
  There is then an exact sequence of $R$--modules,
  \begin{displaymath}
    0 \to M \to H \to A \to 0, 
  \end{displaymath}
  where $A$ is Gorenstein projective and $\pd{H}=\Gpd{M}$.
\end{lem}

\begin{proof}
  If $M$ is Gorenstein projective, we take $0 \to M \to H \to A \to 0$
  to be the first short exact sequence in the ``right half'' of a
  complete projective resolution of $M$.
  
  We may now assume that $\Gpd{M}= n> 0$. By \cite[thm.~2.10]{HH3}
  there exists an exact sequence,
  \begin{displaymath}
    0 \to K \to A' \to M \to 0, 
  \end{displaymath}
  where $A'$ is Gorenstein projective, and $\pd{K}=n-1$. Since $A'$ is
  Gorenstein projective, there exists (as above) a short exact
  sequence,
  \begin{displaymath}
    0 \to A' \to Q \to A \to 0, 
  \end{displaymath}
  where $Q$ is projective, and $A$ is Gorenstein projective. Consider
  the push-out:
  \begin{displaymath} 
    \xymatrix@=1.5em{{} & {} & 0 & 0 & {} \\
      {} & {} & A \ar[u] \ar@{=}[r] & A \ar[u] & {} \\  
      0 \ar[r] & K \ar[r] & Q \ar[u] \ar[r] & H
      \ar[u] \ar[r] & 0 \\
      0 \ar[r] & K \ar@{=}[u] \ar[r] & A' \ar[u] \ar[r] & M
      \ar[u] \ar[r] & 0 \\ 
      {} & {} & 0 \ar[u] & 0 \ar[u] & {} }    
  \end{displaymath}  
  The second column of this diagram is the desired sequence. To see
  this we must argue that $\pd{H} = n$: The class of Gorenstein
  projective modules is projectively resolving by \cite[thm.\ 
  2.5]{HH3}, so if $H$ were projective, exactness of the second column
  would imply that $\Gpd{M}=0$, which is a contradiction.
  Consequently $\pd{H}>0$.  Applying, e.g., \cite[ex.~4.1.2(1)]{Wei}
  to the first row above it follows that $\pd{H} = \pd{K}+1 = n$.
\end{proof}

The next two lemmas have proofs similar to that of
\eqref{lem:helping_sequence_for_GPD}.

\begin{lem}
  \label{lem:helping_sequence_for_GID}
  Let $N$ be an $R$--module of finite Gorenstein injective dimension.
  There is then an exact sequence of $R$--modules,
  \begin{displaymath}
    0 \to B \to H \to N \to 0, 
  \end{displaymath}
  where $B$ is Gorenstein injective and $\id{H}=\Gid{N}$. \qed
\end{lem}

\begin{lem}
  \label{lem:helping_sequence_for_GFD}
  Assume that $R$ is right coherent, and let $M$ be an $R$--module of
  finite Gorenstein flat dimension. There is then an exact sequence of
  $R$--modules,
  \begin{displaymath}
    0 \to M \to H \to A \to 0, 
  \end{displaymath}
  where $A$ is Gorenstein flat and $\fd{H}=\Gfd{M}$. \qed
\end{lem}

%%% SECTION 3
\section{Measuring Gorenstein dimensions}
\label{sec:TTBT}

Gorenstein dimensions are defined in terms of resolutions, and when a
finite resolution is known to exist, the minimal length of such can be
determined by vanishing of certain derived functors. We collect these
descriptions in three theorems, which mimic the style of Cartan and
Eilenberg. Such results have previously in
\cite{LWC3,HH3,EEnOJn93,EEnOJn95,EJT93} been established in more
restrictive settings, and the purpose of this section is to present
them in the more general setting of complexes over associative rings.

We start by investigating the Gorenstein projective dimension.

\begin{thm}
  \label{thm:GPDT}
  Let $X\in\Dbr$ be a complex of finite Gorenstein projective
  dimension. For $n\in\ZZ$ the following are equivalent:
  \begin{eqc}
  \item $\Gpd{X}\le n$.
  \item $n\ge \inf{U} - \inf{\DHom{X}{U}}$ for all $U\in\Db$ of finite
    projective or finite injective dimension with $\H{U} \ne 0$.
  \item $n \ge -\inf{\DHom{X}{Q}}$ for all projective $R$--modules
    $Q$.
  \item $n\ge \sup{X}$ and, for any right-bounded complex $A\eq X$ of
    Gorenstein projective modules, the cokernel $\cC{n}{A} =
    \Coker{(A_{n+1} \lora A_n)}$ is a Gorenstein projective module.
  \end{eqc}
  
  Moreover, the following hold:
  \begin{align*}
    \Gpd{X} &= \supremum{\inf{U}-\inf{\DHom{X}{U}}}{\pd{U}<\infty
      \, \text{ \upshape and } \, \H{U}\ne 0}\\
    &= \supremum{-\inf{\DHom{X}{Q}}}{Q \text{ \upshape is
        projective}}\\
    &\le \FPD + \sup{X}.
  \end{align*}
\end{thm}

\begin{proof}
  The proof of the equivalence of \eqclbl{i}--\eqclbl{iv} is cyclic.
  Clearly, \eqclbl{ii} is stronger than \eqclbl{iii}, and this leaves
  us three implications to prove.
  
  \proofofimp{i}{ii} Choose a complex $A\in\Cb$ consisting of
  Gorenstein projective modules, such that $A\eq X$ and $A_\l =0$ for
  $\l >n$. First, let $U$ be a complex of finite projective dimension
  with $\H{U} \ne 0$. Set $i=\inf{U}$ and note that $i \in \ZZ$ as $U
  \in \Db$ with $\H{U} \ne 0$.  Choose a bounded complex $P\eq U$ of
  projective modules with $P_\l=0$ for $\l<i$. By corollary
  \eqref{prp:rpst-CovDHom} the complex $\Hom{A}{P}$ is isomorphic to
  $\DHom{X}{U}$ in $\D[\ZZ]$; in particular,
  $\inf{\DHom{X}{U}}=\inf{\Hom{A}{P}}$. For $\l<i-n$ and $q\in\ZZ$,
  either $q>n$ or $q+\l\le n+\l<i$, so the module
  \begin{gather*}
    \Hom{A}{P}_\l \qeq \prod_{q\in\ZZ}\Hom{A_q}{P_{q+\l}}
  \end{gather*}
  vanishes. Hence, $\H[\l]{\Hom{A}{P}}=0$ for $\l<i-n$, and
  $\inf{\DHom{X}{U}} \ge i-n = \inf{U}-n$ as desired.
  
  Next, let $U$ be a complex of finite injective dimension and choose
  a bounded complex $I\eq U$ of injective modules. Set $i=\inf{U}$ and
  consider the soft truncation $V = \Tsr{i}{I}$. The modules in $V$
  have finite injective dimension and $U\eq V$, whence $\Hom{A}{V} \eq
  \DHom{X}{U}$ by corollary \eqref{prp:rpst-CovDHom}, and the proof
  continues as above.
  
  \proofofimp{iii}{iv} This part evolves in three steps. First we
  establish the inequality $n \ge \sup{X}$, next we prove that the
  $n$'th cokernel in a bounded complex $A\eq X$ of Gorenstein
  projectives is again Gorenstein projective, and finally we give an
  argument that allows us to conclude the same for $A\in\Cbr$.
  
  To see that $n \ge \sup{X}$, it is sufficient to show that
  \begin{align*}
    \tag{\one} \supremum{-\inf{\DHom{X}{Q}}}{Q \text{ \upshape is
        projective}} \ge \sup{X}.
  \end{align*}
  By assumption, $g=\Gpd{X}$ is finite; i.e., $X\eq A$ for some
  complex
  \begin{align*}
    A = 0\to A_g \xra{\dd{g}{A}} A_{g-1}\to \dots \to A_i \to 0,
  \end{align*}
  and it is clear $g \ge \sup{X}$ since $X\eq A$. For any projective
  module $Q$, the complex $\Hom{A}{Q}$ is concentrated in degrees $-i$
  to $-g$,
  \begin{align*}
    0 \to \Hom{A_i}{Q} \to \dots \to \Hom{A_{g-1}}{Q}
    \xra{\Hom{\dd{g}{A}}{Q}} \Hom{A_{g}}{Q} \to 0,
  \end{align*}
  and isomorphic to $\DHom{X}{Q}$ in $\D[\ZZ]$, cf.~corollary
  \eqref{prp:rpst-CovDHom}. First, consider the case $g=\sup{X}$: The
  differential $\mapdef{\dd{g}{A}}{A_g}{A_{g-1}}$ is not injective, as
  $A$ has homology in degree $g=\sup{X}=\sup{A}$. By the definition of
  Gorenstein projective modules, there exists a projective module $Q$
  and an injective homomorphism $\mapdef{\varphi}{A_g}{Q}$. Because
  $\dd{g}{A}$ is not injective, $\varphi \in \Hom{A_{g}}{Q}$ cannot
  have the form $\Hom{\dd{g}{A}}{Q}(\psi) = \psi\dd{g}{A}$ for some
  $\psi \in \Hom{A_{g-1}}{Q}$. That is, the differential
  $\Hom{\dd{g}{A}}{Q}$ is not surjective; hence $\Hom{A}{Q}$ has
  non-zero homology in degree $-g = -\sup{X}$, and (\one) follows.
  
  Next, assume that $g > \sup{X}= s$ and consider the exact sequence
  \begin{align*}
    \tag{\three} 0 \to A_g \to \dots \to A_{s+1} \to A_s \to \cC{s}{A}
    \to 0.
  \end{align*}
  It shows that $\Gpd{\cC{s}{A}} \le g-s$, and it is easy to check
  that equality must hold; otherwise, we would have $\Gpd{X}<g$.  A
  straightforward computation based on
  corollary~\eqref{prp:rpst-CovDHom}, cf.~\cite[lem.~(4.3.9)]{LWC3},
  shows that for all $m>0$, all $n\ge \sup{X}$, and all projective
  modules $Q$ one has
  \begin{gather*}
    \tag{\two} \Ext{m}{\cC{n}{A}}{Q}=\H[-(m+n)]{\DHom{X}{Q}}.
  \end{gather*}
  By \cite[thm.~2.20]{HH3} we have $\Ext{g-s}{\cC{s}{A}}{Q} \ne 0$ for
  some projective $Q$, whence $\H[-g]{\DHom{X}{Q}} \ne 0$ by (\two),
  and (\one) follows. We conclude that $n \ge \sup{X}$.
  
  It remains to prove that $\cC{n}{A}$ is Gorenstein projective for
  any right-bounded complex $A\eq X$ of Gorenstein projective modules.
  By assumption, $\Gpd{X}$ is finite, so a bounded complex
  $\widetilde{A} \eq X$ of Gorenstein projective modules does exist.
  Consider the cokernel $\cC{n}{\widetilde{A}}$.  Since $n \ge \sup{X}
  = \sup{\widetilde{A}}$, it fits in an exact sequence $0 \to
  \widetilde{A}_t \to \dots \to \widetilde{A}_{n+1} \to
  \widetilde{A}_n \to \cC{n}{\widetilde{A}} \to 0$, where all the
  $\widetilde{A}_\l$'s are Gorenstein projective. By (\two) and
  \cite[thm.~2.20]{HH3} it now follows that also
  $\cC{n}{\widetilde{A}}$ is Gorenstein projective.  With this, it is
  sufficient to prove the following:
  \begin{quote}
    If $P,A\in\Cbr$ are complexes of, respectively, projective and
    Gorenstein projective modules, and $P\eq X\eq A$, then the
    cokernel $\cC{n}{P}$ is Gorenstein projective if and only if
    $\cC{n}{A}$ is so.
  \end{quote}
  Let $A$ and $P$ be two such complexes.  As $P$ consists of
  projectives, there is a quasi-isomorphism
  $\mapdef[\xre]{\pi}{P}{A}$, cf.~\cite[1.4.P]{LLAHBF91}, which
  induces a quasi-isomorphism between the truncated complexes,
  $\mapdef[\xre]{\Tsl{n}{\pi}}{\Tsl{n}{P}}{\Tsl{n}{A}}$.  The mapping
  cone
  \begin{align*}
    \Mc{(\Tsl{n}{\pi})} = 0 \to \cC{n}{P} \to P_{n-1} \oplus \cC{n}{A}
    \to P_{n-2} \oplus A_{n-1} \to \cdots
  \end{align*}
  is a bounded exact complex, in which all modules but the two
  left-most ones are known to be Gorenstein projective modules. It
  follows by the resolving properties of Gorenstein projective
  modules, cf.~\cite[thm.~2.5]{HH3}, that $\cC{n}{P}$ is Gorenstein
  projective if and only if $P_{n-1} \oplus \cC{n}{A}$ is so, which is
  tantamount to $\cC{n}{A}$ being Gorenstein projective.
  
  \proofofimp{iv}{i} Choose a projective resolution $P$ of $X$; by
  \eqclbl{iv} the truncation $\Tsl{n}{P}$ is a complex of the desired
  type.
  
  To show the last claim, we still assume that $\Gpd{X}$ is finite.
  The two equalities are immediate consequences of the equivalence of
  \eqclbl{i}, \eqclbl{ii}, and \eqclbl{iii}. Moreover, it is easy to
  see how a complex $A\in\Cb$ of Gorenstein projective modules, which
  is isomorphic to $X$ in $\D$, may be truncated to form a Gorenstein
  projective resolution of the top homology module of $X$,
  cf.~(\three) above.  Thus, by the definition we automatically obtain
  the inequality $\Gpd{X} \le \FGPD + \sup{X}$, where
  \begin{displaymath}
    \FGPD = \sup \left\{ \, \Gpd{M} \,
      \left|
        \begin{array}{l}
          \mbox{$M$ is an $R$--module with finite} \\
          \mbox{Gorenstein projective dimension}
        \end{array}
      \right.
    \right\}
  \end{displaymath}
  is the (left) finitistic Gorenstein projective dimension, cf.\ 
  \eqref{finitistic dimensions}. Finally, we have $\FGPD = \FPD$ by
  \cite[thm.~2.28]{HH3}.
\end{proof}

\begin{cor}
  \label{cor:GPDT_ring}
  Assume that $R$ is left coherent, and let $X\in\Dbr$ be a complex
  with finitely presented homology modules. If $X$ has finite
  Gorenstein projective dimension, then
  \begin{align*}
    \Gpd{X} \, = \, -\inf{\DHom{X}{R}}.
  \end{align*}
\end{cor}

\begin{proof}
  Under the assumptions, $X$ admits a resolution by finitely generated
  projective modules, say $P$; and thus, $\Hom{P}{-}$ commutes with
  arbitrary sums. The proof is now a straightforward computation.
\end{proof}

Next, we turn to the Gorenstein injective dimension. The proof of
theorem \eqref{thm:GIDT} below relies on
corollary~\eqref{prp:rist-CovDHom} instead of \eqref{prp:rpst-CovDHom}
but is otherwise similar to the proof of theorem \eqref{thm:GPDT};
hence it is omitted.

\begin{thm}
  \label{thm:GIDT}
  Let $Y\in\Dbl$ be a complex of finite Gorenstein injective
  dimension. For $n\in\ZZ$ the following are equivalent:
  \begin{eqc}
  \item $\Gid{Y}\le n$.
  \item $n\ge -\sup{U} - \inf{\DHom{U}{Y}}$ for all $U\in\Db$ of
    finite injective or finite projective dimension with $\H{U}\ne 0$.
  \item $n \ge -\inf{\DHom{J}{Y}}$ for all injective $R$--modules $J$.
  \item $n \ge -\inf{Y}$ and, for any left-bounded complex $B\eq Y$ of
    Gorenstein injective modules, the kernel $\cZ{-n}{B} =
    \Ker{(B_{-n} \lora B_{-(n+1)})}$ is a Gorenstein injective module.
  \end{eqc}
  
  Moreover, the following hold:
  \begin{xxalignat}{3}
    && \Gid{Y} &= \supremum{-\sup{U}-\inf{\DHom{U}{Y}}}{\id{U}<\infty
      \, \text{ \upshape and } \, \H{U}\ne 0}\\
    &&&= \supremum{-\inf{\DHom{J}{Y}}}{J \text{ \upshape is injective}}\\
    &&&\le \FID - \inf{Y}. &&\qed
  \end{xxalignat}
\end{thm}

The next result is a straightforward application of Matlis' structure
theorem for injective modules to the equality in \eqref{thm:GIDT}.

\begin{cor} \label{cor:GIDT_hulls}
  Assume that $R$ is commutative and noetherian. If $\,Y\in\Dbl$ is a
  complex of finite Gorenstein injective dimension, then
  \begin{xxalignat}{3}
    && \Gid{Y} \, &= \, \supremum{\,-\inf{\DHom{\E{R/\p}}{Y}}\,}{ \p
      \in \SpecR}. &&\qed
  \end{xxalignat}
\end{cor}

Finally, we treat the Gorenstein flat dimension.

\begin{thm}
  \label{thm:GFDT}
  Assume that $R$ is right coherent, and let $X\in\Dbr$ be a complex
  of finite Gorenstein flat dimension. For $n\in\ZZ$ the following are
  equivalent:
  \begin{eqc}
  \item $\Gfd{X}\le n$.
  \item $n\ge \supP{\Dtp{U}{X}} - \sup{U}$ for all $U\in\Db[\Ropp]$ of
    finite injective or finite flat dimension with $\H{U}\ne 0$.
  \item $n\ge \supP{\Dtp{J}{X}}$ for all injective $\Ropp$--modules
    $J$.
  \item $n\ge \sup{X}$ and, for any right-bounded complex $A\eq X$ of
    Gorenstein flat modules, the cokernel $\cC{n}{A} = \Coker{(A_{n+1}
      \lora A_n)}$ is Gorenstein flat.
  \end{eqc}
  
  Moreover, the following hold:
  \begin{align*}
    \Gfd{X} &=
    \supremum{\supP{\Dtp{U}{X}}-\sup{U}}{\id[\Ropp]{U}<\infty
      \, \text{ \upshape and } \, \H{U}\ne 0}\\
    &= \supremum{\supP{\Dtp{J}{X}}}{J \text{ \upshape is injective}}\\
    &\le \FFD + \sup{X}.
  \end{align*}
\end{thm}

\begin{proof}
  The proof of the equivalence of \eqclbl{i}--\eqclbl{iv} is cyclic.
  The implication \proofofimp[]{ii}{iii} is immediate, and this leaves
  us three implications to prove.
  
  \proofofimp{i}{ii} Choose a complex $A\in\Cb$ of Gorenstein flat
  modules, such that $A\eq X$ and $A_\l =0$ for $\l >n$. First, let
  $U\in\D[\Ropp]$ be a complex of finite injective dimension with
  $\H{U} \ne 0$. Set $s=\sup{U}$ and pick a bounded complex $I\eq U$
  of injective modules with $I_\l=0$ for $\l>s$. By corollary
  \eqref{prp:rfst-CovDHom} the complex $\tp{I}{A}$ is isomorphic to
  $\Dtp{U}{X}$ in $\D[\ZZ]$; in particular,
  $\supP{\Dtp{U}{X}}=\supP{\tp{I}{A}}$. For $\l>n+s$ and $q\in\ZZ$
  either $q>s$ or $\l-q\ge \l-s>n$, so the module
  \begin{gather*}
    \tpP{I}{A}_\l \qeq \coprod_{q\in\ZZ} \tp{I_q}{A_{\l -q}}
  \end{gather*}
  vanishes. Hence, $\H[\l]{\tp{I}{A}}=0$ for $\l>n+s$, forcing
  $\supP{\Dtp{U}{X}} \le n+s = n+ \sup{U}$ as desired.
  
  Next, let $U\in\D[\Ropp]$ be a complex of finite flat dimension and
  choose a bounded complex $F\eq U$ of flat modules. Set $s=\sup{U}$
  and consider the soft truncation $V = \Tsl{s}{F}$. The modules in
  $V$ have finite flat dimension and $U\eq V$, hence $\tp{V}{A} \eq
  \Dtp{U}{X}$ by corollary \eqref{prp:rfst-CovDHom}, and the proof
  continues as above.
  \enlargethispage*{1cm}
  \proofofimp{iii}{iv} By assumption, $\Gfd{X}$ is finite, so a
  bounded complex $A \eq X$ of Gorenstein flat modules does exist.
  For any injective $\Ropp$--module $J$, we have $\Dtp{J}{X} \eq
  \tp{J}{A}$ in $\D[\ZZ]$ by corollary \eqref{prp:rfst-CovDHom}, so
  \begin{align*}
    \supP{\Dtp{J}{X}} &= \supP{\tp{J}{A}}\\
    &= -\inf{\Hom[\ZZ]{\tp{J}{A}}{\QQ/\ZZ}}\\
    &= -\inf{\Hom[\Ropp]{J}{\Hom[\ZZ]{A}{\QQ/\ZZ}}}\\
    &= -\inf{\DHom[\Ropp]{J}{\Hom[\ZZ]{A}{\QQ/\ZZ}}},
  \end{align*}
  where the last equality follows from corollary
  \eqref{prp:rist-CovDHom}, as $\Hom[\ZZ]{A}{\QQ/\ZZ}$ is a complex of
  Gorenstein injective modules by \cite[thm.~3.6]{HH3}. As desired, we
  now have:
  \begin{align*}
    n &\ge \supremum{\supP{\Dtp{J}{X}}}{J \text{ \upshape is
        injective}}\\
    &=\supremum{-\inf{\DHom[\Ropp]{J}{\Hom[\ZZ]{A}{\QQ/\ZZ}}}}{J
      \text{
        \upshape is injective}}\\
    &\ge -\inf{\Hom[\ZZ]{A}{\QQ/\ZZ}}\\
    &= \sup{A} = \sup{X},
  \end{align*}
  where the inequality follows from \eqref{thm:GIDT} (applied to
  $\Ropp$). The rest of the argument is similar to the one given in
  the proof of theorem \eqref{thm:GPDT}. It uses that the class of
  Gorenstein flat modules is resolving, and here we need the
  assumption that $R$ is right-coherent, cf.~\cite[thm.~3.7]{HH3}.

  \proofofimp{iv}{i} Choose a projective resolution $P$ of $X$; by
  \eqclbl{iv} the truncation $\Tsl{n}{P}$ is a complex of the desired
  type.
  
  For the second part, we can argue, as we did in the proof of
  theorem~\eqref{thm:GPDT}, to see that $\Gfd{X} \le \FGFD + \sup{X}$,
  where
  \begin{displaymath}
    \FGFD = \sup \left\{ \, \Gfd{M} \,
      \left|
        \begin{array}{l}
          \mbox{$M$ is an $R$--module with finite} \\
          \mbox{Gorenstein flat dimension}
        \end{array}
      \right.
    \right\}
  \end{displaymath}
  is the (left) finitistic Gorenstein flat dimension,
  cf.~\eqref{finitistic dimensions}. By~\cite[thm.~3.24]{HH3} we have
  $\FGFD = \FFD$, and this concludes the proof.
\end{proof}

The next corollary is immediate by Matlis' structure theorem for
injective modules.

\begin{cor}
  \label{cor:GFDT_hulls}
  Assume that $R$ is commutative and noetherian. If $X\in\Dbr$ is a
  complex of finite Gorenstein flat dimension, then
  \begin{xxalignat}{3}
    && \Gfd{X} \, &= \, \supremum{\,\supP{\Dtp{\E{R/\p}}{X}}}{ \p \in
      \SpecR}. &&\qed
  \end{xxalignat}
\end{cor}

The next two results deal with relations between the Gorenstein
projective and flat dimensions.  Both are Gorenstein versions of
well-established properties of the classical homological dimensions.

\begin{prp}
  \label{Gfd<Gpd}
  Assume that $R$ is right coherent and any flat $R$--module has
  finite projective dimension. For every $X \in \Dbr$ the next
  inequality holds
  \begin{align*}
    \Gfd{X} \le \Gpd{X}.
  \end{align*}
\end{prp}

\begin{proof}
  Under the assumptions, it follows by \cite[proof of prop.~3.4]{HH3}
  that every Gorenstein projective $R$--module also is Gorenstein
  flat.
\end{proof}

We now compare the Gorenstein projective and Gorenstein flat dimension
to Auslander and Bridger's G--dimension. In \cite{MAuMBr69} Auslander
and Bridger introduce the G--dimension, $\Gdim(-)$, for finitely
generated modules over a left and right noetherian ring.

The G--dimension is defined via resolutions consisting of modules from
the so-called G--class, $\G$. The G--class consists exactly of the
finite $R$--modules $M$ with $\Gdim{M} = 0$ (together with the
zero-module). The basic properties are catalogued in
\cite[prop.~(3.8)(c)]{MAuMBr69}.

When $R$ is commutative and noetherian, \cite[sec.~2.3]{LWC3}
introduces a G--dimension, also denoted $\Gdim(-)$, for complexes in
$\Dfbr$. For modules it agrees with Auslander and Bridger's
G--dimension. However, the definition given in \cite[sec.~2.3]{LWC3}
makes perfect sense over any two-sided noetherian ring.

\begin{prp}
  \label{prp:Gdim}
  Assume that $R$ is left and right coherent. For a complex \mbox{$X\in
  \Dbr$} with finitely presented homology modules, the following hold:
  \begin{prt}
  \item If every flat $R$--module has finite projective dimension, then
    \begin{displaymath}
      \Gpd{X} \, = \, \Gfd{X}.
    \end{displaymath}
  \item If $R$ is left and right noetherian, then
    \begin{displaymath}
      \Gpd{X} \, = \, \Gdim{X}.
    \end{displaymath}
  \end{prt}
\end{prp}

\begin{proof}
  Since $R$ is right coherent and flat $R$--modules have finite
  projective dimension, proposition \eqref{Gfd<Gpd} implies that
  $\Gfd{X} \le \Gpd{X}$. To prove the opposite inequality in
  \prtlbl{a}, we may assume that $n=\Gfd{X}$ is finite. Since $R$ is
  left coherent and the homology modules of $X$ are finitely
  presented, we can pick a projective resolution $P$ of $X$, where
  each $P_\l$ is finitely generated. The cokernel $\cC{n}{P}$ is
  finitely presented, and by theorem~\eqref{thm:GFDT} it is Gorenstein
  flat.
  
  Following the proof of \cite[thm.~(5.1.11)]{LWC3}, which deals with
  commutative, noetherian rings and is propelled by Lazard's
  \cite[lem.~1.1]{DLz69}, it is easy, but tedious, to check that over
  a left coherent ring, any finitely presented Gorenstein flat module
  is also Gorenstein projective. Therefore, $\cC{n}{P}$ is actually
  Gorenstein projective, which shows that $\Gpd{X} \le n$ as desired.
  
  Next, we turn to \prtlbl{b}. By the ``if'' part of
  \cite[thm.~(4.2.6)]{LWC3}, every module in the G--class is
  Gorenstein projective in the sense of definition \eqref{dfn:gproj}.
  (Actually, \cite[thm.~(4.2.6)]{LWC3} is formulated under the
  assumption that $R$ is commutative and noetherian, but the proof
  carries over to two-sided noetherian rings as well.) It follows
  immediately that $\Gpd{X} \le \Gdim{X}$.
  
  For the opposite inequality, we may assume that $n=\Gpd{X}$ is
  finite. Let $P$ be a projective resolution of $X$ by finitely
  generated modules, and consider cokernel the $\cC{n}{P}$.  Of
  course, $\cC{n}{P}$ is finitely generated, and by
  theorem~\eqref{thm:GPDT} it is also Gorenstein projective.  Now the
  ``only if'' part of (the already mentioned ``associative version''
  of) \cite[thm.~(4.2.6)]{LWC3} gives that $\cC{n}{P}$ belongs to the
  G--class. Hence, $\Tsl{n}{P}$ is a resolution of $X$ by modules from
  the G--class and, thus, $\Gdim{X} \leq n$.
\end{proof}

\enlargethispage*{1cm}
\begin{rmk}
  It is natural to ask if finiteness of Gorenstein dimensions is
  ``closed under distinguished triangles''. That is, in a
  distinguished triangle,
  \begin{displaymath}
    X \to Y \to Z \to \Sigma X, 
  \end{displaymath}
  where two of the three complexes $X$, $Y$ and $Z$ have finite, say,
  Gorenstein projective dimension, is then also the third complex of
  finite Gorenstein projective dimension?
  
  Of course, once we have established the main theorems,
  \eqref{thm:mainthm_A} and \eqref{thm:mainthm_B}, it follows that
  over a ring with a dualizing complex, finiteness of each of the
  three Gorenstein dimensions is closed under distinguished triangles.
  This conclusion is immediate, as the Auslander categories are
  triangulated subcategories of $\D$.  However, from the definitions
  and results of this section, it is not immediately clear that the
  Gorenstein dimensions possess this property in general. However, in
  \cite{OV03} Veliche introduces a Gorenstein projective dimension for
  unbounded complexes. By \cite[thm.~3.2.8(1)]{OV03}, finiteness of
  this dimension is closed under distinguished triangles; by
  \cite[thm.~3.3.6]{OV03} it coincides, for right-bounded complexes,
  with the Gorenstein projective dimension studied in this paper.
\end{rmk}

%%% SECTION 4
\section{Auslander categories}
\label{sec:ABC}

In this section, we prove two theorems linking finiteness of
Gorenstein homological dimensions to Auslander categories:

\begin{thm}
  \label{thm:mainthm_A}
  Let $\pair$ be a noetherian pair with a dualizing complex $\bi{D}$.
  For $X \in \Dbr$ the following conditions are equivalent:
  \begin{eqc}
  \item $X \in \A$.
  \item $\Gpd X$ is finite.
  \item $\Gfd X$ is finite.
  \end{eqc}
\end{thm}

\begin{cor}
  \label{cor:GorGRJ}
  Let $R$ be commutative noetherian with a dualizing complex, and let
  $X\in\Dbr$. Then $\Gfd{X}$ is finite if and only if $\,\Gpd{X}$ is
  finite. \qed
\end{cor}

\begin{cor}
  Let $R$ and $S$ be commutative noetherian local rings and
  $\mapdef{\varphi}{R}{S}$ be a local homomorphism. If $R$ has a
  dualizing complex, then $\Gdim[\no]{\varphi}$ is finite if and only
  if $\,\Gfd{S}$ is finite.
\end{cor}

\begin{proof}
  As $R$ admits a dualizing complex, then \cite[thm.~(4.3)]{LLAHBF97}
  yields that $\Gdim[]{\varphi}$ is finite precisely when $S \in \A$.
  It remains to invoke theorem \eqref{thm:mainthm_A}.
\end{proof}

\begin{thm}
  \label{thm:mainthm_B}
  Let $\pair$ be a noetherian pair with a dualizing complex $\bi{D}$.
  For $Y \in \Dbl[S]$ the following conditions are equivalent:
  \begin{eqc}
  \item $Y \in \B[S]$.
  \item $\Gid[S]{Y}$ is finite.
  \end{eqc}
\end{thm}

At least in the case $R=S$, this connection between Auslander
categories and Gorenstein dimensions has been conjectured/expected.
One immediate consequence of the theorems above is that the full
subcategory, of $\D$, of complexes of finite Gorenstein
projective/flat dimension is equivalent, cf.~\eqref{bfhpg:Foxby_eq},
to the full subcategory, of $\D[S]$, of complexes of finite Gorenstein
injective dimension.

The main ingredients of the proofs of theorems \eqref{thm:mainthm_A}
and \eqref{thm:mainthm_B} are lemmas \eqref{lem:suf_Gproj} and
\eqref{lem:suf_Ginj}, respectively. However, we begin with the
following:

\begin{lem}
  \label{lem:FD}
  Let $\pair$ be a noetherian pair with a dualizing complex $\bi{D}$.
  For $X \in \A$ and $Y \in \B[S]$ the following hold:
  \begin{prt}
  \item For all $R$--modules $M$ with finite $\fd{M}$ there is an
    inequality,
    \begin{displaymath}
      -\inf{\DHom{X}{M}} \le \id[S]{\biP{D}} + \supP{\Dtp{\bi{D}}{X}}.
    \end{displaymath}
  \item For all $\Ropp$--modules $M$ with finite $\id[\Ropp]{M}$ there
    is an inequality,
    \begin{displaymath}
      \supP{\Dtp{M}{X}} \le \id[S]{\biP{D}} + \supP{\Dtp{\bi{D}}{X}}.
    \end{displaymath}
  \item For all $S$--modules $N$ with finite $\id[S]{N}$ there is an
    inequality,
    \begin{displaymath}
      -\inf{\DHom[S]{N}{Y}} \le \id[S]{\biP{D}} - \inf{\DHom[S]{\bi{D}}{Y}}.
    \end{displaymath}
  \end{prt}
\end{lem}

\begin{proof}
  \proofoftag{a} If $\H{X}=0$ or $M=0$ there is nothing to prove.
  Otherwise, we compute as follows:
  \begin{align*}
    -\inf{\DHom{X}{M}} &\,=\,
    -\inf{\DHom{X}{\DHom[S]{\bi{D}}{\Dtp{\bi{D}}{M}}}} \\
    &\,=\,
    -\inf{\DHom[S]{\Dtp{\bi{D}}{X}}{\Dtp{\bi{D}}{M}}} \\
    &\,\le\,
    \id[S]{(\Dtp{\bi{D}}{M})} + \supP{\Dtp{\bi{D}}{X}} \\
    &\,\le\, \id[S]{\biP{D}} - \inf{M} + \supP{\Dtp{\bi{D}}{X}}
  \end{align*}
  The first equality follows as $M \in \A$ and the second by
  adjointness. The first inequality is by \cite[thm.~2.4.I]{LLAHBF91};
  the second is by \cite[thm.~4.5(F)]{LLAHBF91}, as $S$ is left
  noetherian and $\fd{M}$ is finite.
  
  \proofoftag{b} Because $\id[\Ropp]{M}$ is finite and $R$ is right
  noetherian, \cite[thm.~4.5(I)]{LLAHBF91} implies that the
  $R$--module $\Hom[\ZZ]{M}{\QQ/\ZZ} \eq \DHom[\ZZ]{M}{\QQ/\ZZ}$ has
  finite flat dimension. Now the desired result follows directly from
  \prtlbl{a}, as
  \begin{align*}
    \supP{\Dtp{M}{X}}
    &\,=\, -\inf \DHom[\ZZ]{\Dtp{M}{X}}{\QQ/\ZZ} \\
    &\,=\, -\inf \DHom{X}{\DHom[\ZZ]{M}{\QQ/\ZZ}}.
  \end{align*}
  
  \proofoftag{c} Again we may assume that $\H{Y} \ne 0$ and $N \ne 0$,
  and hence:
  \begin{align*}
    -\inf{\DHom[S]{N}{Y}} &\,=\,
    -\inf{\DHom[S]{\Dtp{\bi{D}}{\DHom[S]{\bi{D}}{N}}}{Y}} \\
    &\,=\,
    -\inf{\DHom{\DHom[S]{\bi{D}}{N}}{\DHom[S]{\bi{D}}{Y}}} \\
    &\,\le\, \pd{\DHom[S]{\bi{D}}{N}} - \inf{\DHom[S]{\bi{D}}{Y}}.
  \end{align*}
  The first equality follows as $N \in \B[S]$ and the second by
  adjointness.  The inequality follows from
  \cite[thm.~2.4.P]{LLAHBF91}.  Now, since $R$ is right noetherian and
  $\id{N}$ is finite, \cite[thm.~4.5(I)]{LLAHBF91} implies that:
  \begin{displaymath}
    \fd{\DHom[S]{\bi{D}}{N}} \,\le\, \id[\Ropp]{\biP{D}} + \sup{N} \,=\,
    \id[\Ropp]{\biP{D}} < \infty.
  \end{displaymath}
  Therefore proposition \eqref{prp:PJresult} gives the first
  inequality in:
  \begin{align*}
    \pd{\DHom[S]{\bi{D}}{N}} &\,\le\, \max\left\{
      \begin{gathered}
        \id[S]{\biP{D}} + \supP{\Dtp{\bi{D}}{\DHom[S]{\bi{D}}{N}}},\\
        \sup{\DHom[S]{\bi{D}}{N}}
      \end{gathered}
    \right\} \\
    &\,\le\, \max \{\,\id[S]{\biP{D}} + \sup{N}\, , \,
    \sup{N} - \infP{\bi{D}}\} \\
    &\,=\, \id[S]{\biP{D}}.\qedhere
  \end{align*}
\end{proof}

\begin{lem}
  \label{lem:suf_Gproj}
  Let $\pair$ be a noetherian pair with a dualizing complex $\bi{D}$.
  If $M$ is an $R$--module satisfying:
  \begin{prt}
  \item $M \in \A$ and
  \item $\Ext{m}{M}{Q} = 0$ for all integers $m>0$ and all projective
    $R$--modules $Q$,
  \end{prt}
  then $M$ is Gorenstein projective.
\end{lem}

\begin{proof}
  We are required to construct a complete projective resolution of
  $M$. For the left half of this resolution, any ordinary projective
  resolution of $M$ will do, because of \prtlbl{b}.  In order to
  construct the right half, it suffices to construct a short exact
  sequence of $R$--modules,
  \begin{gather}
    \tag{\one} 0 \to M \to P' \to M' \to 0,
  \end{gather}
  where $P'$ is projective and $M'$ satisfies \prtlbl{a} and
  \prtlbl{b}.  The construction of (\one) is done in three steps.
  
  \step{1} \, First we show that $M$ can be embedded in an $R$--module
  of finite flat dimension. Consider resolutions of $\bi{D}$,
  cf.~\eqref{dfn:dc}\prtlbl{2\&3},
  \begin{displaymath}
    \bi{P} \,\xre\, \bi{D} \,\xre\, \bi{I},
  \end{displaymath}
  where $\bi{I}$ is bounded, and let $\lambda \colon \bi{P} \xre
  \bi{I}$ be the composite of these two quasi-isomorphisms.  Since
  $M \in \A$, the complex $\Dtp{\bi{D}}{M} \eq \tp{\bi{P}}{M}$ belongs
  to $\Db[S]$; in particular, $\tp{\bi{P}}{M}$ admits an
  $S$--injective resolution,
  \begin{gather*}
    \tp{\bi{P}}{M} \xre J \in \Cbl[S].
  \end{gather*} 
  Concordantly, we get quasi-isomorphisms of $R$--complexes,
  \begin{gather*}
    \tag{\four} M \xre \Hom[S]{\bi{P}}{\tp{\bi{P}}{M}} \xre
    \Hom[S]{\bi{P}}{J} \xle \Hom[S]{\bi{I}}{J}.
  \end{gather*} 
  Note that since $R$ is right noetherian, and $\bi{I}$ is a complex
  of bimodules consisting of injective $\Ropp$--modules, while $J$ is
  a complex of injective $S$--modules, the modules in the $R$--complex
  \begin{displaymath}
    F = \Hom[S]{\bi{I}}{J} \in \Cbl
  \end{displaymath}
  are flat.  From (\four) it follows that the modules $M$ and
  $\H[0]{F}$ are isomorphic, and that $\H[\l]{F}=0$ for all $\l \neq
  0$.  Now, $\H[0]{F}$ is a submodule of the zeroth cokernel
  $\cC{0}{F}= \Coker{(F_1 \lora F_0)}$, and $\cC{0}{F}$ has finite flat
  dimension over $R$ since
  \begin{gather*}
    \cdots \to F_1 \to F_0 \to \cC{0}{F} \to 0
  \end{gather*}
  is exact and $F \in \Cbl$. This proves the first claim.
  
  \step{2} \, Next, we show that $M$ can be embedded in a flat
  (actually free) $R$--module.  Note that, by induction on $\pd K$,
  condition \prtlbl{b} is equivalent to
  \begin{prt}
  \item[\prtlbl{b'}] $\Ext{m}{M}{K} = 0$ for all $m>0$ and all
    $R$--modules $K$ with $\pd K<\infty$.
  \end{prt}
  By the already established \stepref{\bf{1}} there exists an
  embedding $M \hookrightarrow C$, where $C$ is an $R$--module of
  finite flat dimension.  Pick a short exact sequence of $R$--modules,
  \begin{gather}
    \tag{\two} 0 \to K \to L \to C \to 0,
  \end{gather}
  where $L$ is free and, consequently, $\fd{K}<\infty$.  Proposition
  \eqref{prp:PJresult} implies that also $\pd{K}$ is finite, and hence
  $\Ext{1}{M}{K}=0$ by \prtlbl{b'}. Applying $\Hom{M}{-}$ to (\two),
  we get an exact sequence of abelian groups,
  \begin{gather*} 
    \Hom{M}{L} \to \Hom{M}{C} \to \Ext{1}{M}{K}=0,
  \end{gather*}
  which yields a factorization,
  \begin{gather*} 
    \xymatrix{M \ar[r] \ar@{..>}[dr] & C \\
      {} & L \ar[u]}
  \end{gather*}
  As $M \hookrightarrow C$ is a monomorphism, so is the map from $M$
  into the free $R$--module $L$.
  \enlargethispage*{1cm}
  \step{3} \, Finally, we are able to construct (\one).  Since $R$ is
  right noetherian there exists by \cite[prop.~5.1]{EEn81} a flat
  preenvelope $\varphi \colon M \lora F$ of the $R$--module $M$.  By
  \stepref{\bf{2}}, $M$ can be embedded into a flat $R$--module, and
  this forces $\varphi$ to be a monomorphism.  Now choose a projective
  $R$--module $P'$ surjecting onto $F$, that is,
  \begin{gather*}
    0 \to Z \to P' \xra{\pi} F \to 0
  \end{gather*}
  is exact. Repeating the argument above, we get a factorization
  \begin{gather*} 
    \xymatrix{M \ar[r]^-{\varphi} \ar@{..>}[dr]_-{\partial} & F \\
      {} & P' \ar[u]_-{\pi}}
  \end{gather*}
  and because $\varphi$ is injective so is $\partial$. Thus, we have a
  short exact sequence
  \begin{gather*}
    \tag{\three} 0 \to M \xra{\partial} P' \to M' \to 0.
  \end{gather*}
  What remains to be proved is that $M'$ has the same properties as
  $M$. The projective $R$--module $P'$ belongs to the $\A$, and by
  assumption so does $M$. Since $\A$ is a triangulated subcategory of
  $\D$, also $M'\in\A$.  Let $Q$ be projective; for $m>0$ we have
  $\Ext{m}{M}{Q} = 0 = \Ext{m}{P'}{Q}$, so it follows from the long
  exact sequence of Ext~modules associated to (\three) that
  $\Ext{m}{M'}{Q}=0$ for $m>1$. To see that also $\Ext{1}{M'}{Q}=0$,
  we consider the exact sequence of abelian groups,
  \begin{gather*}
    \Hom{P'}{Q} \xra{\Hom{\partial}{Q}} \Hom{M}{Q} \to \Ext{1}{M'}{Q}
    \to 0.
  \end{gather*}
  Since $Q$ is flat and $\mapdef{\fe}{M}{F}$ is a flat preenvelope,
  there exists, for each\linebreak \mbox{$\tau\in\Hom{M}{Q}$}, a \ho
  $\mapdef{\tau'}{F}{Q}$ such that $\tau=\tau'\fe$; that is,
  $\tau=\tau'\pr\partial = \Hom{\partial}{Q}(\tau'\pr)$.  Thus, the
  induced map $\Hom{\partial}{Q}$ is surjective and, therefore,
  $\Ext{1}{M'}{Q}=0$.
\end{proof}

\begin{proof}[Proof of theorem \eqref{thm:mainthm_A}]
  \proofofimp{ii}{iii} By proposition \eqref{prp:PJresult}, every flat
  $R$--module has finite projective dimension. Furthermore, $R$ is
  right noetherian, and thus $\Gfd{X} \le \Gpd{X}$ by proposition
  \eqref{Gfd<Gpd}.
  
  \proofofimp{iii}{i} If $\Gfd X$ is finite, then, by definition, $X$
  is isomorphic in $\D$ to a bounded complex $A$ of Gorenstein flat
  modules. Consider resolutions of the dualizing complex,
  cf.~\eqref{dfn:dc}\prtlbl{2\&3},
  \begin{displaymath}
    \bi{P} \,\xre\, \bi{D} \,\xre\, \bi{I},
  \end{displaymath}
  where $\bi{I}$ is bounded, and let $\lambda \colon \bi{P} \xre
  \bi{I}$ be the composite quasi-isomorphism.  As
  $\id[\Ropp]{\biP{D}}$ is finite, theorem~\eqref{thm:GFDT} implies
  that $\Dtp{\bi{D}}{X}$ is bounded. Whence, to prove that $X \in \A$,
  we only need to show that
  \begin{gather*}
    \eta_A \colon A \longrightarrow \Hom[S]{\bi{P}}{\tp{\bi{P}}{A}}
  \end{gather*}
  is a quasi-isomorphism. Even though the modules in ${}_SP$ are not
  necessarily finitely generated, we do have ${}_SP \eq {}_SD \in
  \Dfb[S]$ by assumption. Since $S$ is left noetherian, there exists a
  resolution,
  \begin{displaymath}
    \xymatrix{\Cbr[S] \ni L \ar[r]^-{\sigma}_-{\eq} & {}_SP}
  \end{displaymath}
  by finitely generated free $S$--modules. There is a commutative
  diagram, in $\C[\ZZ]$,
  \begin{gather*}
    \xymatrix{ \tp{\Hom[S]{{}_SP}{\bi{I}}}{A}
      \ar[d]_-{\tp{\Hom[S]{\sigma}{\bi{I}}}{A}}^-{\eq}
      \ar[rr]^-{\textnormal{tensor-eval.}} & {} &
      \Hom[S]{{}_SP}{\tp{\bi{I}}{A}}
      \ar[d]_-{\eq}^-{\Hom[S]{\sigma}{\tp{\bi{I}}{A}}} \\
      \tp{\Hom[S]{L}{\bi{I}}}{A}
      \ar[rr]_-{\textnormal{tensor-eval.}}^-{\simeq} & {} &
      \Hom[S]{L}{\tp{\bi{I}}{A}} }
  \end{gather*}
  Since both $L$ and ${}_SP$ are right bounded complexes of projective
  modules, the \qiso $\sigma$ is preserved by the functor
  $\Hom[S]{-}{U}$ for any $S$--complex $U$.  This explains why the
  right vertical map in the diagram above is a \qiso.  Since $L \in
  \Cbr[S]$ consists of finitely generated free $S$--modules and
  $\bi{I}$ and $A$ are bounded, it follows by, e.g.,
  \cite[lem.~4.4.(F)]{LLAHBF91} that the lower horizontal
  tensor-evaluation morphism is an isomorphism in the category of
  $\ZZ$--complexes.  Finally, $\Hom[S]{\sigma}{\bi{I}}$ is a \qiso
  between complexes in $\Cbl[\Ropp]$ consisting of injective modules.
  This can be seen by using the so-called swap-isomorphism:
  \begin{displaymath}
    \Hom[\Ropp]{-_R}{\Hom[S]{{}_SP}{\bi{I}}} \cong
    \Hom[S]{{}_SP}{\Hom[\Ropp]{-_R}{\bi{I}}}.
  \end{displaymath}
  Now theorem \eqref{thm:P->I}\prtlbl{b} implies that also
  $\tp{\Hom[S]{\sigma}{\bi{I}}}{A}$ is a \qiso.  This argument proves
  that the lower horizontal tensor-evaluation map in the next
  commutative diagram of $R$--complexes is a \qiso:
  \begin{gather*}
    \xymatrix{\tp{{}_RR_R}{A} \is  A \ar[rr]^-{\eta_A} \ar[d]_-{\tp{(\Hom[S]{\lambda}{\bi{I}} \circ
          \rhty[I])}{A}}^-{\simeq} & {} &
      \Hom[S]{\bi{P}}{\tp{\bi{P}}{A}}
      \ar[d]_-{\simeq}^-{\Hom[S]{\bi{P}}{\tp{\lambda}{A}}} \\
      \tp{\Hom[S]{\bi{P}}{\bi{I}}}{A}
      \ar[rr]_-{\textnormal{tensor-eval.}}^-{\simeq} & {} &
      \Hom[S]{\bi{P}}{\tp{\bi{I}}{A}} }
  \end{gather*}
  It remains to see that the vertical morphisms in the above diagram
  are invertible:
  \begin{itemlist}
  \item Consider the composite $\gamma$ of the following two
    quasi-isomorphisms of complexes of $(R,\Ropp)$--\-bimodules,
    cf.~the appendix,
    \begin{displaymath}
      \xymatrix{{}_RR_R \ar[r]_-{\eq}^-{\rhty[I]} & \Hom[S]{\bi{I}}{\bi{I}}
        \ar[rr]_-{\eq}^-{\Hom[S]{\lambda}{\bi{I}}} & {} & \Hom[S]{\bi{P}}{\bi{I}}.
      }
    \end{displaymath}
    First note that $R$ and $\Hom[S]{\bi{P}}{\bi{I}}$ belong to
    $\Cbl[\Ropp]$.  Clearly, $R$ is a flat $\Ropp$--module; and we
    have already seen that $\Hom[S]{\bi{P}}{\bi{I}}$ consists of
    injective $\Ropp$--modules.  Therefore, theorem
    \eqref{thm:P->I}\prtlbl{b} implies that $\tp{\gamma}{A}$ is a
    quasi-isomorphism.\footnote{Note that the results in
      section~\ref{sec:Ubiquity} do not allow us to conclude that the
      individual morphisms $\tp{\rhty[I]}{A}$ and
      $\tp{\Hom[S]{\lambda}{\bi{I}}}{A}$ are quasi-isomorphisms.}
    
  \item Since $\bi{P}, \bi{I} \in \Cbr[\Ropp]$, and $\bi{P}$ consists
    of projective $\Ropp$--modules, while $\bi{I}$ consists of
    injective $\Ropp$--modules, it follows by theorem
    \eqref{thm:P->I}\prtlbl{a} that the induced morphism,
    $\tp{\lambda}{A} \colon \tp{\bi{P}}{A} \xre \tp{\bi{I}}{A}$, is a
    \qiso.  Hence also $\Hom[S]{\bi{P}}{\tp{\lambda}{A}}$ is a \qiso.
  \end{itemlist}
  
  \proofofimp{i}{ii} Let $X \in \A$; we can assume that $\H{X}\ne 0$.
  By lemma \eqref{lem:FD},
  \begin{displaymath}
    -\inf{\DHom{X}{M}} \leq s = \id[S]{\biP{D}} + \supP{\Dtp{\bi{D}}{X}}
  \end{displaymath}
  for all $R$--modules $M$ with $\fd{M}<\infty$. Set
  $n=\max\{s,\sup{X}\}$.  Take a projective resolution $\Cbr\ni P\xre
  X$. Since $n \geq \sup{X} = \sup{P}$ we have \mbox{$\Tsl{n}{P} \simeq P
  \simeq X$}, and hence it suffices to show that the cokernel
  $\cC{n}{P} = \Coker{(P_{n+1} \lora P_n)}$ is a Gorenstein projective
  $R$--module.  By lemma~\eqref{lem:suf_Gproj} it is enough to prove
  that
  \begin{prt}
  \item $\cC{n}{P} \in \A$ and
  \item $\Ext{m}{\cC{n}{P}}{Q} = 0$ for all integers $m>0$ and all
    projective $R$--modules $Q$.
  \end{prt}
  Consider the exact sequence of complexes
  \begin{equation*}
    0 \to \Thl{n-1}{P} \to \Tsl{n}{P} \to \Sigma^n \cC{n}{P} \to 0.
  \end{equation*}
  Obviously, $\Thl{n-1}{P}$ belongs to $\A$,
  cf.~\eqref{bfhpg:Foxby_eq}, and also $\Tsl{n}{P} \simeq P \simeq X
  \in \A$.  Because $\A$ is a triangulated subcategory of $\D$, we
  conclude that $\Sigma^n\cC{n}{P}$, and hence $\cC{n}{P}$, belongs to
  $\A$. This establishes \prtlbl{a}.
  
  To verify \prtlbl{b}, we let $m>0$ be an integer, and $Q$ be any
  projective $R$--module. Since $n \ge \sup{X}$, it is a
  straightforward computation, cf.~\cite[lem.~(4.3.9)]{LWC3}, to see
  that
  \begin{gather*}
    \Ext{m}{\cC{n}{P}}{Q} \,\is\, \H[-(m+n)]{\DHom{X}{Q}},
  \end{gather*}
  for $m>0$.  Since $-\inf{\DHom{X}{Q}} \leq s \leq n$, we see that
  $\Ext{m}{\cC{n}{P}}{Q}=0$.
\end{proof}  

\begin{proof}[Proof of Theorem \eqref{thm:mainthm_B}]
  Using lemma \eqref{lem:suf_Ginj} below, the proof is similar to that
  of theorem \eqref{thm:mainthm_A}.  Just as the proof of lemma
  \eqref{lem:suf_Gproj} uses $R$--flat preenvelopes, the proof of
  lemma \eqref{lem:suf_Ginj} below uses $S$--injective precovers.  The
  existence of such precovers is guaranteed by \cite{MLT76},
  cf.~\cite[prop.~2.2]{EEn81}, as $S$ is left noetherian.
\end{proof}

\begin{lem}
  \label{lem:suf_Ginj}
  Let $\pair$ be a noetherian pair with a dualizing complex $\bi{D}$.
  If $N$ is an $S$--module satisfying
  \begin{prt}
  \item $N \in \B[S]$ and
  \item $\Ext[S]{m}{J}{N} = 0$ for all integers $m>0$ and all
    injective $S$--modules $J$,
  \end{prt}
  then $N$ is Gorenstein injective.\qed
\end{lem}

%%% SECTION 5
\section{Stability results}
\label{sec:Stability}

We now apply the characterization from the previous section to show
that finiteness of Gorenstein dimensions is preserved under a series
of standard operations. In this section, all rings are commutative and
noetherian.

\begin{ipg}
  It is known from \cite{LWC3} that $\Gid{\Hom{X}{E}} \le \Gfd{X}$ for
  $X\in\Dbr$ and injective modules $E$. Here is a dual result, albeit
  in a more restrictive setting:
\end{ipg}

\begin{prp}
  \label{prp:duality}
  Let $R$ be commutative and noetherian with a dualizing complex, and
  let $E$ be an injective $R$--module. For $Y \in \Dbl$ there is an
  inequality,
  \begin{displaymath}
    \Gfd \Hom{Y}{E} \leq \Gid{Y}, 
  \end{displaymath}
  and equality holds if $\,E$ is faithful.
\end{prp}

\begin{proof}
  It is sufficient to prove that if $N$ is a Gorenstein injective
  module, then $\Hom{N}{E}$ is Gorenstein flat, and that the converse
  holds, when $E$ is faithful.
  
  Write $-^{\vee} = \Hom{-}{E}$ for short, and set $d = \FFD$, which
  is finite as $R$ has a dualizing complex, cf.~\eqref{finitistic
    dimensions}. From theorem \eqref{thm:GFDT} we know that if $C$ is
  any module with $\Gfd{C} < \infty$ then, in fact, $\Gfd{C} \leq d$.
  
  Now assume that $N$ is Gorenstein injective, and consider part of
  the left half of a complete injective resolution of $N$,
  \begin{equation}
    \tag{\one}
    0 \to C_d \to I_{d-1} \to \cdots \to I_0 \to N \to 0. 
  \end{equation}
  The $I_\l$'s are injective $R$--modules and $C_d$ is Gorenstein
  injective. In particular, \mbox{$C_d \in \B$} by theorem
  \eqref{thm:mainthm_B}, and $C_d^{\,\vee}\in \A$ by
  \cite[lem.~(3.2.9)(b)]{LWC3}, so \mbox{$\Gfd{C_d^{\,\vee}}\le d$}.
  Applying the functor $-^{\vee}$ to (\one) we obtain an exact
  sequence:
  \begin{equation*}
    0 \to N^{\,\vee}\to I_0^{\,\vee} \to \cdots \to I_{d-1}^{\,\vee} \to
    C_d^{\,\vee} \to 0, 
  \end{equation*} 
  where the $I_\l^{\,\vee}$'s are flat $R$--modules. From theorem
  \eqref{thm:GFDT} we conclude that $N^{\,\vee}$ is Gorenstein flat.
  
  Finally, we assume that $E$ is faithfully injective and that
  $N^{\vee}$ is Gorenstein flat, in particular, $N^{\vee} \in \A$.
  This forces $N \in \B$, again by \cite[lem.~(3.2.9)(b)]{LWC3}, that
  is, $\Gid{N}$ is finite. By lemma
  \eqref{lem:helping_sequence_for_GID} there exists an exact sequence
  \begin{equation*}
    0 \to B \to H \to N \to 0,
  \end{equation*}
  where $B$ is Gorenstein injective and $\id{H}=\Gid{N}$. By the first
  part of the proof, $B^{\vee}$ is Gorenstein flat, and by assumption,
  so is $N^{\vee}$.  Therefore, exactness of
  \begin{displaymath}
    0 \to N^{\vee} \to H^{\vee} \to B^{\vee} \to 0
  \end{displaymath}
  forces $H^{\vee}$ to be Gorenstein flat by the resolving property of
  Gorenstein flat modules, cf.~\cite[thm.~3.7]{HH3}. In particular,
  $\Gid{N} = \id{H} = \fd{H^{\vee}} =\Gfd{H^{\vee}}=0$. Here the
  second equality follows from \cite[thm.~1.5]{Ishikawa}.
\end{proof}

The next result is an immediate corollary of \cite[thm.~(6.4.2) and
(6.4.3)]{LWC3} and \eqref{prp:duality}.

\begin{cor}
  \label{cor:tensor}
  Let $R$ be commutative and noetherian with a dualizing complex, and
  let $F$ be a flat $R$--module. For $Y \in \Dbl$ there is an
  inequality,
  \begin{displaymath}
    \Gid \tpP{Y}{F} \leq \Gid{Y}, 
  \end{displaymath}
  and equality holds if $\,F$ is faithful.\qed
\end{cor}

\begin{thm}
  \label{thm:base change}
  Let $\mapdef{\varphi}{R}{S}$ be a homomorphism of commutative
  noetherian rings with $\fd[\no]{\varphi}$ finite. Assume that $R$
  has a dualizing complex $D$ and $E=\Dtp{D}{S}$ is dualizing for $S$.
  For $Y\in\Db$ the following hold:
  \begin{align*}
    \tag{a} \Gfd{Y} < \infty & \Ra \Gfd[S]{(\Dtp{Y}{S})} <
    \infty \\
    \tag{b} \Gid{Y} < \infty & \Ra \Gid[S]{(\Dtp{Y}{S})} < \infty
  \end{align*}
  Either implication may be reversed under each of the next two extra
  conditions:
  \begin{itemlist}
  \item $\varphi$ is faithfully flat.
  \item $\varphi$ is local and the complex $Y$ belongs to $\Dfb$.
  \end{itemlist}
  {\rm When $\mapdef{\varphi}{\Rm}{\Sn}$ is a local homomorphism, the
  assumption that the base-changed complex $\Dtp{D}{S}$ is dualizing
  for $S$ is tantamount to $\varphi$ being Gorenstein (at the maximal
  ideal $\n$ of $S$). For details see \cite[thm.~(7.8)]{LLAHBF97}.}
\end{thm}

\begin{proof}
  We only prove the statements for the Gorenstein injective dimension;
  the proof for the Gorenstein flat dimension is similar.
  
  In view of theorem \eqref{thm:mainthm_B}, we need to see that the
  base changed complex $\Dtp{Y}{S}$ belongs to $\B[S]$ when $Y\in\B$.
  This is a special case of \cite[prop.~(5.9)]{LWC1}, from where it
  also follows that the implication may be reversed when $\varphi$ is
  faithfully flat.
  
  Next, let $\varphi$ be local, $Y$ be in $\Dfb$, and assume that
  $\Dtp{Y}{S} \in \B[S]$. The aim is to show that $Y \in \B$. First,
  we verify that $\DHom{D}{Y}$ has bounded homology. As $E =
  \Dtp{D}{S}$ is a dualizing complex for $S$, we may compute as
  follows
  \begin{align*}
    \tag{\one}
    \begin{split}
      \DHom[S]{E}{\Dtp{Y}{S}} & \eq \DHom[S]{\Dtp{D}{S}}{\Dtp{Y}{S}} \\
      & \eq \DHom{D}{\Dtp{Y}{S}} \\
      & \eq \Dtp{\DHom{D}{Y}}{S}.
  \end{split}
  \end{align*}
  Here the first isomorphism is trivial, the second is adjointness,
  and the third follows from \cite[(A.4.23)]{LWC3}.
  
  The remainder of the proof is built up around two applications of
  Iversen's amplitude inequality, which is now available for unbounded
  complexes \cite[thm.~3.1]{SriHBF}. The amplitude inequality yields
  \begin{align*}
    \tag{\two} \ampP{\DHom{D}{Y}} \le \ampP{\Dtp{\DHom{D}{Y}}{S}},
  \end{align*}
  as $\varphi$ is assumed to be of finite flat dimension.  Here the
  amplitude of a complex $X$ is defined as $\amp{X} = \sup{X} -
  \inf{X}$. From $(\one)$ we read off that the homology of
  $\Dtp{\DHom{D}{Y}}{S}$ is bounded, and by $(\two)$ this shows that
  the homology of $\DHom{D}{Y}$ is bounded as well.
  
  Finally, consider the commutative diagram
  \begin{gather*}
    \xymatrix{\Dtp[S]{E}{\DHom[S]{E}{\Dtp{Y}{S}}}
      \ar[d]_-{\varepsilon_{Y \otimes^{\mathbf{L}}S}}^-{\eq} & {} &
      \Dtp{D}{\DHom{D}{\Dtp{Y}{S}}} \ar[ll]^-{\eq}_-{\gamma_{Y
          \otimes^{\mathbf{L}} S}} \\
      \Dtp{Y}{S} & {} & \Dtp{(\Dtp{D}{\DHom{D}{Y}})}{S} \ar[u]^-{\simeq}_-{D
        \otimes^{\mathbf{L}}\omega_{DYS}}
      \ar[ll]^-{\varepsilon_{Y}\otimes^{\mathbf{L}} S} , }
  \end{gather*}
  where $\gamma_{Y \otimes^{\mathbf{L}}_R S}$ is a natural isomorphism
  induced by adjointness and commutativity of the derived tensor
  product. The diagram shows that $\varepsilon_Y
  \otimes^{\mathbf{L}}_{R} S$ is an isomorphism. As
  $\Dtp{D}{\DHom{D}{Y}}$ has degreewise finitely generated homology,
  we may apply \cite[prop.~2.10]{SriSean} to conclude that
  $\varepsilon_Y$ is an isomorphism as well.
\end{proof}

\begin{bfhpg}[Localization]
  Working directly with the definition of Gorenstein flat modules (see
  \cite[lem.~(5.1.3)]{LWC3}), it is easily verified that the
  inequality
  \begin{equation*}
    \Gfd[R_{\p}]{X_{\p}} \leq \Gfd{X}
  \end{equation*} 
  holds for all complexes $X\in\Dbr$ and all prime ideals $\p$ of $R$.
  
  Turning to the Gorenstein projective and injective dimensions, it is
  natural to ask if they do not grow under localization. When $R$ is
  local and Cohen--Macaulay with a dualizing module, Foxby settled the
  question affirmatively in \cite[cor.~(3.5)]{HBF94}. More recently,
  Foxby extended the result for Gorenstein projective dimension to
  commutative noetherian rings of finite Krull dimension; see
  \cite[(5.5)(b)]{LWCHH}. Unfortunately, it is not clear how to apply
  the ideas of that proof to the Gorenstein injective dimension, but
  there is a partial result:
\end{bfhpg}

\begin{prp}
  \label{prp:Gid_localization}
  Let $R$ be commutative and noetherian with a dualizing complex. For
  any complex $Y \in \Dbl$ and any prime ideal $\p$ of $R$, there is
  an inequality
  \begin{equation*}
    \Gid[R_{\p}]{Y_{\p}} \leq \Gid{Y}.
  \end{equation*}
\end{prp}

\begin{proof}
  It suffices to show that if $N$ is a Gorenstein injective
  $R$--module, then $N_{\p}$ is Gorenstein injective over $R_{\p}$.
  This is proved in the exact same manner as in
  \cite[cor.~(3.5)]{HBF94} using theorem~\eqref{thm:mainthm_B}.
\end{proof}

It is immediate from definition \eqref{dfn:gflat} that a direct sum of
Gorenstein flat modules is Gorenstein flat. It has also been proved
\cite{EEnJAL02} that, over a right coherent ring, a colimit of
Gorenstein flat modules indexed by a filtered set is Gorenstein flat.

One may suspect that also a product of Gorenstein flat modules is
Gorenstein flat. In the sequel this is proved for commutative
noetherian rings with a dualizing complex. To this end, the next lemma
records an important observation.

\begin{lem}
  \label{lem:ABcl}
  Let $R$ be a commutative noetherian ring with a dualizing complex.
  The Auslander categories $\A$ and $\B$ are closed under direct
  products of modules and under colimits of modules indexed by a
  filtered set.
\end{lem}

\begin{proof}
  There are four parts to the lemma; they have similar proofs, so we
  shall only prove the first claim, that $\A$ is closed under set
  indexed products of modules.
  
  Let $D$ be a dualizing complex for $R$, and let $L \xre D$ be a
  resolution of $D$ by finitely generated free modules. Consider a
  family of modules $\{M_i\}_{i \in I}$ from $\A$ and set \mbox{$M =
    \prod_{i \in I}M_i$}. The canonical chain map
  \begin{displaymath}
    \xymatrix{\tp{L}{M} = \tp{L}{(\textstyle{\prod}_{i\in I}M_i)}
      \ar[r]^-{\alpha} & \textstyle{\prod}_{i\in
        I}\tpP{L}{M_i}}  
  \end{displaymath}
  is an isomorphism. This is a straightforward verification; it hinges
  on the fact that the module functors \mbox{$\tp{L_n}{-}$} commute
  with arbitrary products, as the modules $L_n$ are finitely generated
  and free.
  
  For each $i\in I$, the complex $\tp{L}{M_i}$ represents
  $\Dtp{D}{M_i}$, and by \cite[prop.~(4.8)(a)]{LWC1} there are
  inequalities \mbox{$\supP{\tp{L}{M_i}} \le \sup{D}$}. Thus, the
  isomorphism $\alpha$ and the fact that homology commutes with
  products yields
  \begin{displaymath}
    \supP{\Dtp{D}{M}} = \supP{\tp{L}{M}} \le \sup{D}.
  \end{displaymath}
  In particular, $\Dtp{D}{M}$ is in $\Db$.  Next, consider the commutative
  diagram,
  \begin{displaymath}
    \xymatrix{*+<8pt>{\mspace{52mu}\textstyle{\prod}_{i \in I}M_i=M}
      \ar[r]^-{\eta_M} \ar[d]_-{\prod_{i \in I}\eta_{M_i}}^-{\eq} &
      \Hom{L}{\tp{L}{M}} \ar[d]^-{\Hom{L}{\alpha}}_-{\is}\\
      \textstyle{\prod}_{i \in I}\Hom{L}{\tp{L}{M_i}} \ar[r]^-{\is}_-{\beta}&
      \Hom{L}{\textstyle{\prod}_{i \in I}\tpP{L}{M_i}}
    }
  \end{displaymath}
  The canonical map $\beta$ is an isomorphism of complexes, as
  $\Hom{L}{-}$ commutes with products. The map $\textstyle{\prod}_{i
    \in I} \eta_{M_i}$ is a quasi-isomorphism, because each
  $\eta_{M_i}$ is one. The upshot is that $\eta_M$ is a
  quasi-isomorphism, and $M$ belongs to $\A$.
\end{proof}
\enlargethispage*{1cm}
\begin{thm}
  \label{thm:prod_Gflat}
  Let $R$ be commutative and noetherian with a dualizing complex. A
  direct product of Gorenstein flat modules is Gorenstein flat.
\end{thm}

\begin{proof}
  Let $A^{(i)}$ be a family of Gorenstein flat modules. By lemma
  \eqref{lem:ABcl} the product $\prod_i A^{(i)}$ is in $\A$ and,
  therefore, $\Gfd{\prod_i A^{(i)}}$ is finite, in fact, at most
  $d=\FFD$, cf.~theorem \eqref{thm:GFDT}. For each $A^{(i)}$ take a
  piece of a complete flat resolution:
  \begin{equation*}
    0 \to A^{(i)} \to F^{(i)}_0 \to \dots \to F^{(i)}_{1-d} \to
    Z^{(i)}_{-d} \to 0,
  \end{equation*}
  where the $F$'s are flat and the $Z$'s are Gorenstein flat. Taking
  products we get an exact sequence:
  \begin{equation*}
    0 \to \textstyle{\prod_i} A^{(i)} \to \textstyle{\prod_i}F^{(i)}_0
    \to \dots \to \textstyle{\prod_i}F^{(i)}_{1-d} \to
    \textstyle{\prod_i}Z^{(i)}_{-d} \to 0.
  \end{equation*}
  Since $R$ is noetherian, the modules $\prod_i F^{(i)}_{\l}$ are
  flat. As noted above $\prod_i Z^{(i)}_{-d}$ has Gorenstein flat
  dimension at most $d$, which forces the product $\prod_i A^{(i)}$ to
  be Gorenstein flat, cf.~theorem \eqref{thm:GFDT}.
\end{proof}

On a parallel note, it is immediate from definition \eqref{dfn:gim}
that a product of Gorenstein injective modules is Gorenstein
injective. We remark that via theorem \eqref{thm:mainthm_A} this gives
a different proof that $\B$ is closed under direct products of
modules. This shows that information flows in both directions between
Auslander categories and Gorenstein dimensions.

Over a ring with a dualizing complex, the proof above is
easily modified to show that a direct sum of Gorenstein injectives is,
again, Gorenstein injective. In view of lemma \eqref{lem:ABcl} it is
natural to expect that even a colimit of Gorenstein injective modules
will be Gorenstein injective. This is proved in the next section; see
theorem \eqref{thm:Gid_and_limits}.

\begin{bfhpg}[Local (co)homology]
  \label{bfhpg:LCH}
  Let $R$ be commutative and noetherian, and let $\mathfrak{a}$ be an
  ideal of $R$. The right derived local cohomology functor with
  support in $\mathfrak{a}$ is denoted $\RG{(-)}$. Its right adjoint,
  $\LL{(-)}$, is the left derived local homology functor with support
  in $\mathfrak{a}$.  Derived local (co)homology is represented on
  $\D$ as
  \begin{gather*}
    \label{equ:repRG}
    \RG{(-)} \ \eq \ \Dtp{\RG{R}}{-}\ \eq \ 
    \Dtp{\Ca}{-}\quad\text{and} \\
    \LL{(-)} \ \eq \ \DHom{\RG{R}}{-}\ \eq \ \DHom{\Ca}{-},
  \end{gather*}
  where $\Ca$ is the so-called {\v C}ech, or stable Koszul, complex on
  $\mathfrak{a}$; it is defined as follows: Let $a \in R$; the complex
  concentrated in homological degrees $0$ and $-1$:
  \begin{gather*}
    \xymatrix{ \mathrm{C}(a) \; = \; 0 \ar[r] & R \ar[r]^-{\rho_a} &
      R_a \ar[r] & 0, }
  \end{gather*}
  where $R_a$ is the localization of $R$ with respect to $\{ a^n \}_{n
    \geq 0}$ and $\rho_a$ is the natural homomorphism $r \mapsto r/1$,
  is the {\v C}ech complex on $a$. When the ideal $\mathfrak{a}$ is
  generated by $a_1,\ldots,a_n$, the {\v C}ech complex on
  $\mathfrak{a}$ is the tensor product $\bigotimes_{i =
    1}^{n}\mathrm{C}(a_i)$.  Observe that the flat dimension of $\Ca$
  is finite.
  
  The above representations of derived local (co)homology will be used
  without mention in the proofs of theorems \eqref{thm:presv} and
  \eqref{thm:width}.  For local cohomology this representation goes
  back to Grothendieck \cite[prop.~(1.4.1)]{egaIII}; see also
  \cite[lem.~(3.1.1)]{ALL} and the corrections in
  \cite[sec.~1]{PSc03}. Local homology was introduced by Matlis
  \cite[\S~4]{EMt74}, when $\mathfrak{a}$ is generated by a regular
  sequence, and for modules over local Cohen--Macaulay rings the
  representation above is implicit in \cite[thm.~5.7]{EMt74}. The
  general version above is due to Greenlees and May
  \cite[sec.~2]{GreenMay}; see also \cite[sec.~1]{PSc03} for
  corrections.
  
  Since $\Ca$ has finite flat dimension, it is immediate that
  $\RG{(-)}$ preserves homological boundedness as well as finite flat
  and finite injective dimension, see also \cite[thm.~6.5]{HBF79}.
  However, $\Ca$ has even finite projective dimension. This calls for
  an argument: Let $a \in R$ and consider the short exact sequence
  \begin{gather*}
    \xymatrix{ 0 \ar[r] & R[X] \ar[rr]^-{ aX - 1 } && R[X]
      \ar[r]^-{\alpha} & R_a \ar[r] & 0, }
  \end{gather*}
  where $\alpha$ is the homomorphism $f(X) \mapsto f(1/a)$. This short
  exact sequence is a bounded free resolution of $R_a$, whence the
  projective dimension of $R_a$ is at most one.  Let $L_a$ be the
  complex
  \begin{gather*}
    \xymatrix@R=0.0ex@M=0.2ex{    {}  &  R[X]   &{}&{}   &               & {}\\
      L_a \; = \; 0 \ar[r] & \oplus \ar[rrr]^-{
        \left(\begin{array}{cc}\!aX -1 & \iota\!
                                                 \end{array}\right)}
                                             &         &{}&{}      R[X] \ar[r]    & 0.\\
                                             {} & R &{}&{} & {} & {} }
  \end{gather*}
  concentrated in homological degrees $0$ and $-1$, where $\iota$
  denotes the natural embedding of $R$ into $R[X]$.  It is
  straightforward to verify that $L_a$ is a bounded free resolution of
  the {\v C}ech complex $\mathrm{C}(a)$. Thus, if the ideal
  $\mathfrak{a}$ is generated by $a_1,\ldots,a_n$, then $L =
  \bigotimes_{i = 1}^{n}L_{a_i}$ is a bounded free resolution of
  $\Ca$.  This shows that the projective dimension of $\Ca$ is finite.
\end{bfhpg}
  
  The last two results investigate preservation of finite Gorenstein
  dimensions by local (co)homology functors.
%\enlargethispage*{1cm}
\begin{thm}
  \label{thm:presv}
  Let $R$ be a commutative noetherian ring, and let $\mathfrak{a}$ be
  an ideal of $R$. For $Y \in \Db$ the following hold:
  \begin{align*}
    \tag{a}
    \label{equ:RGGfd}
    \Gfd{Y}<\infty  \Ra \Gfd{(\RG{Y})}\, < \, \infty.
  \end{align*}
  If, in addition, $R$ admits a dualizing complex, then
  \begin{align*}
    \tag{b}
    \label{equ:RGGid}
    \Gid{Y}<\infty  \Ra \Gid{(\RG{Y})}\, <\, \infty.
  \end{align*}
  Moreover, if $R$ has a dualizing complex, both implications may be
  reversed if $\mathfrak{a}$ is in the Jacobson radical of $R$ and $Y
  \in \Dfb$.
\end{thm}

\begin{ipg}
  Note that in theorem \eqref{thm:presv} we use the existence of a
  dualizing complex to establish preservation of finite Gorenstein
  injective dimension. In \eqref{thm:presv_LL} below the dualizing
  complex is used to establish preservation of finite Gorenstein flat
  dimension.
\end{ipg}
\enlargethispage*{1cm}
\begin{thm}
  \label{thm:presv_LL}
  Let $R$ be a commutative noetherian ring, and let $\mathfrak{a}$ be
  an ideal of $R$. For $X \in \Db$ the following hold:
  \begin{align*}
    \tag{a}
    \label{equ:LLGid}
    \Gid{X}<\infty  \Ra \Gid{(\LL{X})}\, < \, \infty.
  \end{align*}
  If, in addition, $R$ admits a dualizing complex, then
  \begin{align*}
    \tag{b}
    \label{equ:LLGfd}
    \Gfd{X}<\infty \Ra \Gfd{(\LL{X})}\, <\, \infty.
  \end{align*}
  Moreover, both implications may be reversed if $\Rm$ is local and
  complete in its $\m$--adic topology, $\mathfrak{a}$ is in the
  Jacobson radical of $R$, and $X \in
  \mathsf{D}^{\mathrm{art}}_{\sqsubset\mspace{-13mu}\sqsupset}(R)$;
  that is, $X$ has bounded homology and its individual homology
  modules are artinian.
\end{thm}

\begin{proof}[Proof of theorem \eqref{thm:presv}]
  Since $\RG{Y} \eq \Dtp{\Ca}{Y}$ the implication \prtlbl{a} follows
  from \cite[thm.~(6.4.5)]{LWC3}, and \prtlbl{b} from a routine
  application of corollary \eqref{cor:tensor}.
  
  Now assume that $R$ has a dualizing complex $D$ and $\mathfrak{a}$
  is in the Jacobson radical. The arguments showing that the
  implications in \prtlbl{a} and \prtlbl{b} can be reversed are
  similar; we only write out the details for the latter.
  
  Assume that $Y \in \Dfb$ and $\RG{Y}\in\B$. First, we show that
  $\DHom{D}{Y}$ is bounded.  By \eqref{bfhpg:LCH} the projective
  dimension of $\Ca$ is finite and, therefore, $\LL{(-)} \eq
  \DHom{\Ca}{-}$ preserves homological boundedness. In particular,
  $\LL{\DHom{D}{\RG{Y}}}$ is bounded since, already,
  $\DHom{D}{\RG{Y}}$ is. Observe that
  \begin{align*}
    \LL{\DHom{D}{\RG{Y}}} &\eq \DHom{D}{\LL{\RG{Y}}} \\
    &\eq \DHom{D}{\LL{Y}} \\
    &\eq \DHom{D}{\Dtp{Y}{R^{\,\widehat{}}_{\mathfrak{a}}}} \\
    &\eq \Dtp{\DHom{D}{Y}}{R^{\,\widehat{}}_{\mathfrak{a}}}.
  \end{align*}
  Here the first isomorphism is swap, the second is by
  \cite[cor.~(5.1.1)(i)]{ALL}, the third follows from
  \cite[prop.~(2.7)]{F}, and the last is by
  \cite[lem.~4.4(F)]{LLAHBF91}.  As $\mathfrak{a}$ is in the Jacobson
  radical of $R$, the completion $R^{\,\widehat{}}_{\mathfrak{a}}$ is
  a faithful flat $R$--module by \cite[thm.~8.14]{Mat} and, therefore,
  $\DHom{D}{Y}$ is bounded.
  
  To show that $\mapdef{\varepsilon_Y}{\Dtp{D}{\DHom{D}{Y}}}{Y}$ is
  invertible, we consider the commutative diagram
  \begin{gather*}
    \xymatrix{\LL{\RG{Y}} \ar[dd]_-{\beth_Y}^-{\simeq} & {} &
      \LL{(\Dtp{D}{\DHom{D}{\RG{Y}}})} \ar[ll]_-{
        \operatorname{\mathbf{L}\Lambda^{\mathfrak{a}}\varepsilon_{
            \mathbf{R}\Gamma_{\mathfrak{a}}Y}}}^-{\eq}\\
      {} \ar[d] & {} & \LL{\RG{(\Dtp{D}{\DHom{D}{Y}})}}
      \ar[u]^-{\simeq}
      \ar[d]_-{\simeq}^-{\beth_{D\otimes^{\mathrm{\mathbf{L}}}{\mathrm{\mathbf{R}Hom}}(D,Y)}}
      \\
      \LL{Y} & {} & \LL{(\Dtp{D}{\DHom{D}{Y}})}
      \ar[ll]^-{\operatorname{\mathbf{L}\Lambda^{\mathfrak{a}}\varepsilon}_Y},
    }
  \end{gather*}
  The top horizontal morphism is invertible as $\RG{Y}\in\B$. The
  vertical morphisms $\beth$ are invertible, again by
  \cite[cor.~(5.1.1)(i)]{ALL}, and the third vertical morphism is
  induced by tensor evaluation, cf.~\cite[lem.~4.4(F)]{LLAHBF91}. The
  diagram shows that $\LL{\varepsilon_Y}$ is an isomorphism. Now,
  since $\Dtp{D}{\DHom{D}{Y}}$ belongs to $\Dfbr$, cf.~\cite[(1.2.1)
  and (1.2.2)]{LLAHBF97}, it follows by \cite[prop.~(2.7)]{F} that we
  may identify $\LL{\varepsilon_Y}$ with
  $\tp{\varepsilon_Y}{R^{\,\widehat{}}_{\mathfrak{a}}}$. Whence,
  $\varepsilon_Y$ is an isomorphism by faithful flatness of
  $R^{\,\widehat{}}_{\mathfrak{a}}$.
\end{proof}

\begin{proof}[Proof of theorem \eqref{thm:presv_LL}]
  Assume that the Gorenstein injective dimension of $X$ is finite. Let
  $L$ be the bounded free resolution of $\Ca$ described in
  \eqref{bfhpg:LCH}. By assumption there exists a bounded complex, say
  $B$, consisting of Gorenstein injective modules and quasi-isomorphic
  to $X$.  We may represent $\LL{X}$ by the bounded complex
  $\Hom{L}{B}$. It is readily seen that the individual modules in the
  latter complex consist of products of Gorenstein injective modules.
  Consequently, they are Gorenstein injective themselves; see
  \cite[thm.~2.6]{HH3}. In particular, $\LL{X}$ has finite Gorenstein
  injective dimension.
  
  In the presence of a dualizing complex, a similar argument applies
  when the Gorenstein flat dimension of $X$ is finite.  This time we
  use that a product of Gorenstein flat modules is Gorenstein flat by
  \eqref{thm:prod_Gflat}.
  
  As in the proof of theorem \eqref{thm:presv} we only argue why the
  implication in \prtlbl{b} can be reversed; the arguments for
  reversing the implication in \prtlbl{a} are similar.
  
  When $\Rmk$ is complete in its $\m$--adic topology it admits a
  dualizing complex \cite[V.10.4]{RAD}. Moreover, Matlis duality
  \cite[thm.~18.6(v)] {Mat} and the assumption that $X$ has bounded
  artinian homology yields $X \eq X^{\vee \vee}$ where $-^{\vee} =
  \Hom{-}{\E{k}}$ is the Matlis duality functor. Here, $\E{k}$ is the
  injective hull of the residue field $k$. By \cite[(2.10)]{F} we have
  \begin{equation*}
    \tag{\one}
    \LL{X} \eq {\RG{(X^{\vee})}}^{\vee}.
  \end{equation*}
  As the complex $X^{\vee}$ has finite homology, see
  \cite[thm.~18.6(v)]{Mat}, and the functor $-^{\vee}$ is faithful
  and exact, we have the following string of biimplications
  \begin{align*}
    \LL{X} \in \B & \Lra \RG{(X^{\vee})} \in \A \\
    & \Lra X^{\vee} \in \A \\
    & \Lra X \in \B.
  \end{align*}
  Here the first biimplication follows from (\one) in conjunction with
  \cite[lem.~(3.2.9)(a)]{LWC3}; the second follows from theorem
  \eqref{thm:presv} and the third from \cite[lem.~(3.2.9)(a)]{LWC3}.
\end{proof}

\begin{obs}
  Let $R$ be a commutative noetherian ring with a dualizing
  complex, and let $X \in \Db$. We will demonstrate that $\LL{X}$ has
  finite, say, Gorenstein injective dimension, when and only when
  $\RG{X}$ has finite Gorenstein injective dimension. The argument is
  propelled by the isomorphisms
  \begin{equation*}
    \tag{\two}
    \RG{X} \xre \RG{\LL{X}}  \quad\text{and}\quad  \LL{\RG{X}} \xre \LL{X},
  \end{equation*}
  which are valid for any $X \in \D$; for details consult
  \cite[cor.~(5.1.1)]{ALL}. The next string of implications
  \begin{equation*}
    \LL{X} \in \B \Ra \RG{\LL{X}} \in \B \Ra 
    \RG{X} \in \B,
  \end{equation*}
  where the first follows from \eqref{thm:presv} and the second from
  $(\two)$, together with
  \begin{equation*}
    \RG{X} \in \B \Ra \LL{\RG{X}} \in \B \Ra 
    \LL{X} \in \B,
  \end{equation*}
  where the first follows from \eqref{thm:presv_LL} and the second
  from $(\two)$, prove the claim. A similar argument is available
  for Gorenstein flat dimension.
\end{obs}

%%% SECTION 6
\section{Bass and Chouinard formulas}
\label{sec:Bass}

The theorems in section~\ref{sec:TTBT} give formulas for measuring
Gorenstein dimensions. We close the paper by establishing alternative
formulas that allow us to measure or even compute Gorenstein injective
dimension. In this section $R$ is a commutative noetherian ring.

\begin{bfhpg}[Width]
  \label{dfn:width}
  Recall that when $\Rmk$ is a local ring, the \emph{width} of an
  $R$--complex $X\in\Dbr$ is defined as $\wdt{X} = \infP{\Dtp{k}{X}}$.
  There is always an inequality,
  \begin{equation}
    \label{equ:width-inf}
    \wdt{X} \ge \inf{X},
  \end{equation}
  and by Nakayama's lemma, equality holds for $X \in \Dfbr$.
\end{bfhpg}

\begin{obs}
  \label{obs:ABdagger}
  Let $R$ be a commutative and noetherian ring with a dualizing
  complex. It is easy to see that the functor $-^\dagger$,
  cf.~\eqref{bfhpg:dagger duality}, maps $\A$ to $\B$ and vice versa,
  cf.~\cite[lem.~(3.2.9)]{LWC3}. Consequently,
  \begin{align*}
    \Gid{Y}<\infty \Lra \Gfd{Y^{\dagger}}<\infty \Lra 
    \Gpd{Y^{\dagger}}<\infty
  \end{align*}
  for $Y\in\Dfb$ by theorems~\eqref{thm:mainthm_A} and
  \eqref{thm:mainthm_B} and duality \eqref{equ:dagger duality}.
\end{obs}

\begin{thm}
  \label{thm:bass}
  Let $R$ be a commutative noetherian local ring. For $Y \in \Db$
  there is an inequality,
  \begin{equation*}
    \Gid{Y} \geq \dptR - \wdt{Y}.
  \end{equation*}
  If, in addition, $R$ admits a dualizing complex, and $Y\in\Dfb$ is a
  complex of finite Gorenstein injective dimension, then the next equality
  holds,
  \begin{equation*}
    \Gid{Y} = \dptR - \inf{Y}.
  \end{equation*}
  In particular,
  \begin{align*}
    \Gid{N} = \dptR
  \end{align*}
  for any finitely generated $R$--module $N \neq 0$ of finite
  Gorenstein injective dimension.
\end{thm}

\begin{proof}
  Set $d = \dptR$ and pick an $R$--regular sequence $\defx[d]$. Note
  that the module $T = R/(\x)$ has $\pd{T} = d$. We may assume that
  $\Gid{Y}<\infty$, and the desired inequality now follows from the
  computation:
  \begin{align*}
    \Gid{Y} &\ge -\inf{\DHom{T}{Y}} \\
    &\ge -\wdt{\DHom{T}{Y}} \\
    &= \pd{T}- \wdt{Y}\\
    &=\dptR - \wdt{Y}.
  \end{align*}
  The first inequality follows from theorem \eqref{thm:GIDT} and the
  second from \eqref{equ:width-inf}. The first equality is by
  \cite[thm.~(4.14)(b)]{RHD} and the last by definition of $T$.
  
  Next, assume that $R$ admits a dualizing complex, and that
  $Y\in\Dfb$ with $\Gid{Y}$ finite. It suffices to prove the
  inequality
  \begin{equation*}
    \Gid{Y} \leq \dptR - \inf{Y}.
  \end{equation*}
  
  Duality \eqref{equ:dagger duality} yields $\Gid{Y} =
  \Gid{Y^{\dagger \, \dagger}}$, and by theorem \eqref{thm:GIDT} there
  exists an injective $R$--module $J$, such that $\Gid{Y^{\dagger \,
      \dagger}} = -\inf{\DHom{J}{Y^{\dagger \, \dagger}}}$. In the
  computation,
  \begin{align*}
    \Gid{Y} &= -\inf{\DHom{Y^{\dagger}}{J^{\dagger}}} \\
    &\le \Gpd{Y^{\dagger}} - \inf{J^{\dagger}} \\
    &\le \Gpd{Y^{\dagger}} + \id{D},
  \end{align*}
  the first (in)equality is by swap and the second by theorem
  \eqref{thm:GPDT}, as $J^{\dagger}$ is a complex of finite flat
  dimension and hence of finite projective dimension. The final
  inequality is by \cite[thm.~2.4.I]{LLAHBF91}. Recall from
  \eqref{obs:ABdagger} that $\Gpd{Y^{\dagger}}$ is finite; since
  $Y^{\dagger}$ has bounded and degreewise finitely generated
  homology, it follows from proposition \eqref{prp:Gdim}(b) and
  \cite[thm.~(2.3.13)]{LWC3} that $\Gpd{Y^{\dagger}} = \dptR -
  \dpt{Y^{\dagger}}$. Thus, we may continue as follows:
  \begin{align*}
    \Gid{Y} &\leq \dptR - \dpt{Y^{\dagger}} + \id{D}\\
    &= \dptR - \inf{Y} - \dpt{D} + \id{D}\\
    &= \dptR - \inf{Y}.
  \end{align*}
  Both equalities stem from well-known properties of dualizing
  complexes, see \cite[(3.1)(a) and (3.5)]{LWC1}.
\end{proof}

A dualizing complex $D$ for a commutative noetherian local ring $R$
is said to be \emph{normalized} if $\inf{D} = \dptR$, see
\cite[(2.5)]{LLAHBF97}. There is an equality $\dpt{Y} =
\inf{Y^{\dagger}}$ for all $Y\in\Dfb$, when the dual $Y^{\dagger}$ is
taken with respect to a normalized dualizing complex, see
\cite[(3.1)(a) and (3.2)(a)]{LWC1}. Any dualizing complex can be
normalized by an appropriate suspension.

\begin{cor}
  \label{cor:Gdim_dagger}
  Let $R$ be a commutative noetherian and local ring, and let $D$ be a
  normalized dualizing complex for $R$. The next equalities hold for
  $Y \in \Dfb$,
  \begin{equation*}
    \Gid{Y} =  \Gpd{Y^{\dagger}} = \Gfd{Y^{\dagger}},
  \end{equation*}
  where $Y^{\dagger} = \DHom{Y}{D}$.
\end{cor}

\begin{proof}
  By observation \eqref{obs:ABdagger} the three dimensions $\Gid{Y}$,
  $\Gpd{Y^{\dagger}}$, and $\Gfd{Y^{\dagger}}$ are simultaneously
  finite, and in this case \eqref{equ:dagger duality} and theorem
  \eqref{thm:bass} give:
  \begin{align*}
    \Gid{Y} = \Gid{Y^{\dagger \, \dagger}} = \dptR - \inf{Y^{\dagger
        \, \dagger}} = \dptR - \dpt{Y^{\dagger}},
  \end{align*}
  where the last equality uses that the dualizing complex is
  normalized. By \cite[thm.~(2.3.13)]{LWC3} we have
  $\Gdim{Y^{\dagger}} = \dptR - \dpt{Y^{\dagger}}$, and
  $\Gdim{Y^{\dagger}} = \Gpd{Y^{\dagger}} = \Gfd{Y^{\dagger}}$ by
  proposition~\eqref{prp:Gdim}.
\end{proof}

Theorem \eqref{thm:bass} is a Gorenstein version of Bass' formula for
injective dimension of finitely generated modules. In \cite{LCh76}
Chouinard proved a related formula:
\begin{align*}
  \id{N} \, = \, \supremum{\,\dptRp - \wdt[R_\p]{N_\p}}{\p \in \SpecR}
\end{align*}
for any module $N$ of finite injective dimension over a commutative
noetherian ring. In the following, we extend this formula to
Gorenstein injective dimension over a ring with dualizing complex.  The
first result in this direction is inspired by
\cite[thm.~8.6]{SriSean}.

\begin{thm}
  \label{thm:width}
  Let $\Rmk$ be a commutative noetherian local ring with a
  dualizing complex. Denote by $\E{k}$ the injective hull of the
  residue field.  For a complex $Y\in\Db$ of finite Gorenstein
  injective dimension, the next equality holds,
  \begin{equation*}
    \wdt{Y} = \dptR + \inf{\DHom{\E{k}}{Y}}.
  \end{equation*}
  In particular, $\wdt{Y}$ and $\inf{\DHom{\E{k}}{Y}}$ are
  simultaneously finite.
\end{thm}

\begin{proof}
  By theorem \eqref{thm:mainthm_B}, $Y$ is in $\B$; in particular, $Y
  \eq \Dtp{D}{\DHom{D}{Y}}$. Furthermore, we can assume that $D$ is a
  normalized dualizing complex, in which case
  we have $\RG[\mathfrak{m}]{D} \eq \E{k}$ by \cite[prop.~V.6.1]{RAD}.
  We compute as follows:
  \begin{align*}
    \wdt{Y} &= \wdt{(\Dtp{D}{\DHom{D}{Y}})} \\
    &= \wdt{D} + \wdt{\DHom{D}{Y}} \\
    &= \inf{D} + \inf{\LL[\mathfrak{m}]{\DHom{D}{Y}}} \\
    &= \dptR   + \inf{\DHom{\RG[\mathfrak{m}]{D}}{Y}} \\
    &= \dptR + \inf{\DHom{\E{k}}{Y}}.
  \end{align*}
  The second equality is by \cite[thm.~2.4(b)]{SYs98}, the third is by
  \eqref{equ:width-inf} and \cite[thm.~(2.11)]{F}, while the
  penultimate one is by adjointness, cf.~\eqref{bfhpg:LCH}.
\end{proof}

\begin{cor}
  \label{cor:width}
  Let $\Rmk$ be a commutative noetherian and local ring with a
  dualizing complex. If $N$ is a Gorenstein injective module, then
  \begin{equation*}
    \wdt{N} \ge \dptR,
  \end{equation*}
  and equality holds if $\,\wdt{N}$ is finite.
\end{cor}

\begin{proof}
  We may assume that $\wdt{N}$ is finite. By theorem \eqref{thm:width}
  the module $\Hom{\E{k}}{N} \eq \DHom{\E{k}}{N}$ is non-zero and the
  equality is immediate by the same theorem.
\end{proof}

\begin{cor}
  Let $\Rmk$ be a commutative noetherian and local ring, and let $D$
  be a normalized dualizing complex for $R$. If $\,Y \in \Dfb$ has
  finite Gorenstein injective dimension, then
  \begin{equation*}
    \Gid{Y} \, = \, -\inf{\DHom{\E{k}}{Y}} \,  = \, -\inf{\DHom{D}{Y}}.
  \end{equation*}
\end{cor}

\begin{proof}
  The first equality comes from the computation,
  \begin{displaymath}
    \Gid{Y} \, = \, \dptR - \inf{Y} 
    \, = \, \dptR - \wdt{Y} 
    \, = \, -\inf{\DHom{\E{k}}{Y}},
  \end{displaymath}
  which uses theorem \eqref{thm:bass}, \eqref{equ:width-inf}, and
  theorem \eqref{thm:width}. For the second equality in the corollary,
  we note that
  \begin{align*}
    \DHom{\E{k}}{Y} &\eq \DHom{\RG[\m]{D}}{Y} \\
    &\eq \LL[\m]{\DHom{D}{Y}} \\
    &\eq \Dtp{\DHom{D}{Y}}{\widehat{R}}.
  \end{align*}
  Here the second isomorphism is by adjointness,
  cf.~\eqref{bfhpg:LCH}, and the last one is by \cite[prop.~(2.7)]{F}
  as $\DHom{D}{Y}$ is in $\Dfbr$. Since $\widehat{R}$ is faithfully
  flat, the complexes $\DHom{\E{k}}{Y}$ and $\DHom{D}{Y}$ must have
  the same infimum.
\end{proof}

\begin{thm}
  \label{thm:Chouinard_injective_case}
  Let $R$ be a commutative and noetherian ring with a dualizing
  complex. If $\,Y\in\Db$ has finite Gorenstein injective dimension,
  then
  \begin{displaymath} 
    \Gid{Y} \, = \, \supremum{\,\dptRp - \wdt[R_\p]{Y_\p}}{\p \in
      \SpecR}. 
  \end{displaymath}
\end{thm}

\begin{proof}
  First we show ``$\ge$''. For any prime ideal $\p$ of $R$,
  proposition \eqref{prp:Gid_localization} and theorem
  \eqref{thm:bass} give the desired inequality,
  \begin{displaymath}
    \Gid{Y} \, \ge \, \Gid[R_\p]{Y_\p} \, \ge \, \dptRp -
    \wdt[R_\p]{Y_\p}.
  \end{displaymath}
  
  For the converse inequality, ``$\le$'', we may assume that $\H{Y}
  \neq 0$. Set $s=\sup{Y}$ and $g=\Gid{Y}$; by theorem
  \eqref{thm:GIDT} we may assume that $Y$ has the form
  \begin{equation*}
    0 \to I_s \to I_{s-1} \to \dots \to I_{-g+1} \to B_{-g} \to 0,
  \end{equation*}
  where the $I$'s are injective and $B_{-g}$ is Gorenstein injective.
  Proving the inequality amounts to finding a prime ideal $\p$ of $R$
  such that $\wdt[R_{\p}]{Y_\p} \le \dptRp - g$.
  
  \textbf{Case $\pmb{s=-g}$}. For any integer $m$,
  \begin{equation*}
    \tag{\one}
    \Gid{(\Sigma^m Y)} \, = \, \Gid{Y}-m \quad \text{ and } \quad
    \wdt[R_\p]{(\Sigma^m Y)_\p} \, = \, \wdt[R_\p]{Y_\p}+m,
  \end{equation*}
  so we can assume that $s=-g=0$, in which case $Y$ is a Gorenstein
  injective module. By \cite[lem.~2.6]{HBF79} there exists a prime
  ideal $\p$, such that the homology of $\Dtp[R_\p]{k(\p)}{Y_\p}$ is
  non-trivial; in particular,
  \begin{equation*}
    \wdt[R_\p]{Y_\p} = \inf{(\Dtp[R_\p]{k(\p)}{Y_\p})}
  \end{equation*}
  is finite. By prop.~\eqref{prp:Gid_localization} the $R_\p$--module
  $Y_\p$ is Gorenstein injective, so $\wdt[R_{\p}]{Y_\p} = \dptRp$ by
  corollary \eqref{cor:width}.
  
  \textbf{Case $\pmb{s=-g+1}$}. We may, cf.~(\one), assume that $s=-1$
  and $g=2$. That is,
  \begin{equation*}
    Y = 0 \to I_{-1} \xra{\b} B_{-2} \to 0.
  \end{equation*}
  From the complete injective resolution of the Gorenstein injective
  $B_{-2}$, we get a short exact sequence
  \begin{equation*}
    0 \to B' \to H_{-2} \xra{\a} B_{-2} \to 0,
  \end{equation*}
  where $H_{-2}$ is injective and $B'$ is Gorenstein injective.
  Consider the pull-back:
  \begin{equation*}
    \xymatrix{0 \ar[r] & B' \ar@{=}[d] \ar[r] & B_{-1} \ar[d]^{\b'}
      \ar[r]^{\a'} & I_{-1} \ar[d]^{\b} \ar[r] & 0 \\
      0 \ar[r] & B' \ar[r] & H_{-2} \ar[r]^{\a} & B_{-2} \ar[r] & 0 }
  \end{equation*}
  The rows are exact sequences, and the top one shows that the module
  $B_{-1}$ is Gorenstein injective. By the snake lemma $\b'$ and $\b$
  have isomorphic kernels and cokernels; that is, $\Ker{\b'} \is
  \Ker{\b} = \H[-1]{Y}$ and $\Coker{\b'} \is \Coker{\b} = \H[-2]{Y}$.
  Thus, the homomorphisms $\a'$ and $\a$ make up a quasi-isomorphism
  of complexes, and we have $Y \eq 0 \to B_{-1} \xra{\b'} H_{-2} \to
  0$.  Similarly, there is a short exact sequence of modules $0 \to B
  \to H_{-1} \xra{\g} B_{-1} \to 0$, where $H_{-1}$ is injective and
  $B$ is Gorenstein injective. The diagram
  \begin{equation*}
    \xymatrix{0 \ar[r] & B \ar[d] \ar[r] & H_{-1} \ar[d]^{\g} \ar[r]^{\b'\g}
      & H_{-2} \ar@{=}[d] \ar[r] & 0 \\
      & 0 \ar[r] & B_{-1} \ar[r]^{\b'} & H_{-2} \ar[r] & 0 }
  \end{equation*}
  is commutative, as $\Ker{\g} = B$, and shows that the vertical maps
  form a surjective morphism of complexes. The kernel of this morphism
  ($0 \to B \to B \to 0$) is exact, so it is a homology isomorphism,
  and we may replace $Y$ with the top row of the diagram. Set $H = 0
  \to H_{-1} \to H_{-2} \to 0$, then the natural inclusion of $H$ into
  $Y$ yields a short exact sequence of complexes $0 \to H \to Y \to B
  \to 0$. By construction $\id{H} \le 2$. To see that equality holds,
  let $J$ be an $R$--module such that $\H[-2]{\DHom{J}{Y}} \ne 0$,
  cf.~theorem \eqref{thm:GIDT}, and inspect the exact sequence of
  homology modules
  \begin{equation*}
    \H[-2]{\DHom{J}{H}} \to \H[-2]{\DHom{J}{Y}} \to
    \H[-2]{\DHom{J}{B}} = 0.
  \end{equation*}
  By the Chouinard formula for injective dimension
  \cite[thm.~2.10]{SYs98} we can now choose a prime ideal $\p$ such
  that $\dptRp - \wdt[R_\p]{H_\p} = 2$.  Set $d = \dptRp$ and consider
  the exact sequence:
  \begin{equation*}
    \H[d-1]{\Dtp[R_\p]{k(\p)}{B_\p}} \to
    \H[d-2]{\Dtp[R_\p]{k(\p)}{H_\p}} \to \H[d-2]{\Dtp[R_\p]{k(\p)}{Y_\p}}.
  \end{equation*}
  The left-hand module is 0 by \eqref{prp:Gid_localization} and
  \eqref{cor:width}, while the middle one is non-zero by choice of
  $\p$. This forces $\H[d-2]{\Dtp[R_\p]{k(\p)}{Y_\p}}\ne 0$ and
  therefore $\wdt[R_\p]{Y_\p} \le d-2$ as desired.
  
  \textbf{Case $\pmb{s>-g+1}$}. We may assume that $g=1$ and $s>0$;
  i.e., $Y$ has the form
  \begin{equation*}
    0 \to I_{s} \to \dots \to I_{1} \to I_{0} \to B_{-1} \to 0,
  \end{equation*}
  where the $I$'s are injective and $B_{-1}$ is Gorenstein injective.
  Set $I = \Thr{1}{Y}$ and $Y'=\Thl{0}{Y}$, then we have an exact
  sequence of complexes $0 \to Y' \to Y \to I \to 0$. Recycling the
  argument applied to $H$ above, we see that $\Gid{Y'} =1$. As $-1 \le
  \sup{Y'} \le 0$, it follows by the preceding cases, that we can
  choose a prime ideal $\p$ such that $\dptRp - \wdt[R_\p]{Y'_\p} =
  1$.  Set $d = \dptRp$ and consider the exact sequence of homology
  modules
  \begin{equation*}
    \H[d]{\Dtp[R_\p]{k(\p)}{I_\p}} \to
    \H[d-1]{\Dtp[R_\p]{k(\p)}{Y'_\p}} \to \H[d-1]{\Dtp[R_\p]{k(\p)}{Y_\p}}.
  \end{equation*}
  By construction $\id[R_\p]{I_\p} \le -1$, so the left-hand module is
  0 by the classical Chouinard formula. The module in the middle is
  non-zero by choice of $\p$, and this forces
  $\H[d-1]{\Dtp[R_\p]{k(\p)}{Y_\p}}\ne 0$, which again implies
  $\wdt[R_{\p}]{Y_\p} \le d - 1$ as desired.
\end{proof} 

The final result is parallel to theorem \eqref{thm:prod_Gflat}.

\begin{thm} \label{thm:Gid_and_limits}
  Let $R$ be a commutative and noetherian ring with a dualizing
  complex. A colimit of Gorenstein injective $R$--modules indexed by a
  filtered set is Gorenstein injective.
\end{thm}

\begin{proof*}
  Let $B_i \to B_j$ be a filtered, direct system of Gorenstein
  injective modules. By theorem~\eqref{thm:mainthm_B} all the $B_i$'s
  belong to $\B$, so by lemma~\eqref{lem:ABcl} also the colimit
  $\varinjlim B_i$ is in $\B$. That is, $\Gid{\varinjlim B_i} <
  \infty$.  Since tensor products \cite[thm.~A1.]{Mat} and homology
  commute with filtered colimits \cite[thm.~2.6.15]{Wei}, we have
  \begin{displaymath}
    \wdt[R_\p]{(\varinjlim B_i)_\p} \geq 
    \underset{i}{\inf} \{\wdt[R_\p]{(B_i)_\p}\},
  \end{displaymath}  
  for each prime ideal $\p$. By theorem
  \eqref{thm:Chouinard_injective_case} we now have
  \begin{align*}
    \Gid{\varinjlim B_i} &= \supremum{\,\dptRp -
      \wdt[R_\p]{(\varinjlim B_i)_\p}}{\p \in \SpecR} \\
    &\leq \supremum{\,\dptRp -
      \underset{i}{\inf}\{\wdt[R_\p]{(B_i)_\p}\}}{\p \in
      \SpecR} \\
    &= \underset{i}{\sup}\{ \supremum{\,\dptRp -
      \wdt[R_\p]{(B_i)_\p}}{\p \in \SpecR}\} \\
    &= \underset{i}{\sup}\{\Gid{B_i}\} =0.\qedhere
  \end{align*}
\end{proof*}  

The Chouinard formula, theorem~\eqref{thm:Chouinard_injective_case},
plays a crucial role in the proof above. Indeed, it is not clear from
the formulas in \eqref{thm:GIDT} and \eqref{cor:GIDT_hulls} that
$\Gid{\varinjlim B_i} \le \sup_i\{\Gid{B_i}\}$.  In \cite{LCh76}
Chouinard proved a similar formula for modules of finite flat
dimension. Also this has been extended to Gorenstein flat dimension:
for modules in \cite[thm.~3.19]{HH3} and \cite[thm.~(2.4)(b)]{RHD},
and for complexes in \cite[thm.~8.8]{SriSean}.

\appendix
\section*{Appendix: Dualizing complexes}
\label{sec:dc}
\stepcounter{section}

Dualizing complexes over non-commutative rings are a delicate matter.
The literature contains a number of different, although related,
extensions of Grothendieck's original definition \cite[V.\S 2]{RAD} to
the non-commutative realm. Yekutieli \cite{Yekutieli:DCncGA}
introduced dualizing complexes for associative $\ZZ$--graded algebras
over a field. Later, Yekutieli--Zhang \cite{Yekutieli-Zhang:RADC} gave
a definition for pairs of non-commutative algebras over a field which
has been used by, among others, J{\o}rgensen \cite{PJ:Preprint} and
Wu--Zhang \cite{Wu-Zhang:DCncLR}.  Related definitions can be found in
Frankild--Iyengar--J{\o}rgensen \cite{AF-SI-PJ:DDGMGDGA} and Miyachi
\cite{Miyachi:DCMDT}.

\begin{ipg}
  Definition \eqref{dfn:dc} is inspired by Miyachi
  \cite[p.~156]{Miyachi:DCMDT} and constitutes an extension of
  Yekutieli--Zhang's \cite[def.~1.1]{Yekutieli-Zhang:RADC}: They
  consider a noetherian pair $\pair$ of algebras over a field $k$; a
  complex $D \in \Db[{\tp[k]{S}{\Ropp}}]$ is said to be dualizing for
  $\pair$ if:
  \begin{rqm}
  \item[\prtlbl{i}] $D$ has finite injective dimension over $S$ and
    $\Ropp$.
  \item[\prtlbl{ii}] The homology of $D$ is degreewise finitely
    generated over $S$ and $\Ropp$.
  \item[\prtlbl{iii}] The homothety morphisms $S \longrightarrow
    \DHom[\Ropp]{D}{D}$ in $\D[{\tp[k]{S}{\Sopp}}]$ and $R
    \longrightarrow \DHom[S]{D}{D}$ in $\D[{\tp[k]{R}{\Ropp}}]$ are
    isomorphisms.
  \end{rqm} 
  As also noted in \cite{Yekutieli-Zhang:RADC}, condition \prtlbl{i}
  is equivalent to:
  \begin{rqm}
  \item[\prtlbl{i'}] There exists a quasi-isomorphism $D \xre I$ in
    $\Db[{\tp[k]{S}{\Ropp}}]$ such that each $I_\ell$ is injective
    over $S$ and $\Ropp$.
  \end{rqm} 
  Even more is true: The canonical ring homomorphisms $S
  \longrightarrow \tp[k]{S}{\Ropp} \longleftarrow \Ropp$ give
  restriction functors,
  \begin{displaymath}
    \C[S] \longleftarrow \C[{\tp[k]{S}{\Ropp}}]
    \longrightarrow \C[\Ropp].
  \end{displaymath}
  Since $k$ is a field, these restriction functors are exact,
  cf.~\cite[p.~45]{Yekutieli:DCncGA}, and thus they send
  quasi-isomorphisms to quasi-isomorphisms. They also send
  projective/injective modules to projective/injective modules,
  cf.~\cite[lem.~2.1]{Yekutieli:DCncGA}. Consequently, a
  projective/injective resolution of $D$ in $\C[{\tp[k]{S}{\Ropp}}]$
  restricts to a projective/injective resolution of $D$ in $\C[S]$ and
  in $\C[\Ropp]$.  Thus, in the setting of
  \cite{Yekutieli-Zhang:RADC}, there automatically exist
  quasi-isomorphisms of complexes of $(S,\Ropp)$--bimodules,
  \begin{displaymath}
    \Cbr[{\tp[k]{S}{\Ropp}}] \ni P \xre D 
    \quad \text{ and } \quad 
    D \xre  I \in \Cb[{\tp[k]{S}{\Ropp}}] 
  \end{displaymath} 
  such that each $P_\ell$ (respectively, $I_\ell$) is projective
  (respectively, injective) over $S$ and over $\Ropp$. It is also by
  virtue of these biresolutions of $D$ that the morphisms from
  \prtlbl{iii} above make sense, cf.~\cite[p.~52]{Yekutieli:DCncGA}.
  
  In this paper, we work with a noetherian pair of rings, not just a
  noetherian pair of algebras over a field. Without the underlying
  field $k$, the existence of appropriate biresolutions does not come
  for free. Therefore, the existence of such resolutions has been made
  part of the very definition of a dualizing complex,
  cf.~\eqref{dfn:dc}\prtlbl{2\&3}.
\end{ipg}

Also a few remarks about definition \eqref{dfn:dc}(4) are in
order.\footnote{These remarks are parallel to Yekutieli's
  considerations \cite[p.~52]{Yekutieli:DCncGA} about well-definedness
  of derived functors between derived categories of categories of
  bimodules.}  Let $\bi{P} \xre \bi{D}$ and $\bi{D} \xre \bi{I}$ be as
in \eqref{dfn:dc}\prtlbl{2\&3} and let $\lambda \colon \bi{P} \xre
\bi{I}$ be the composite; note that $\lambda$ is $S$-- and
$\Ropp$--linear. Consider the diagram of complexes of
$(S,\Sopp)$--bimodules,
\begin{displaymath}
  \xymatrix{{}_SS_S \ar[d]_-{\lhty[I]} \ar[rrr]^{\lhty[P]} & {} & {} &
    \Hom[\Ropp]{\bi{P}}{\bi{P}} \ar[d]_-{\eq}^{\Hom[\Ropp]{\bi{P}}{\lambda}} \\ 
    \Hom[\Ropp]{\bi{I}}{\bi{I}} \ar[rrr]^-{\eq}_-{\Hom[\Ropp]{\lambda}{\bi{I}}} & {}
    & {} & \Hom[\Ropp]{\bi{P}}{\bi{I}} 
  }
\end{displaymath}
Note that the following facts:
\begin{itemlist}
\item $\lhty[P]$ and $\lhty[I]$ are both $S$--linear and
  $\Sopp$--linear,
\item $\Hom[\Ropp]{\bi{P}}{\lambda}$ is $\Sopp$--linear, and
\item $\Hom[\Ropp]{\lambda}{\bi{I}}$ is $S$--linear
\end{itemlist} 
are immediate consequences of the $(S,\Sopp)$--bistructures on
\begin{displaymath}
  \Hom[\Ropp]{\bi{P}}{\bi{P}} \ \ , \ \ \Hom[\Ropp]{\bi{P}}{\bi{I}}
  \ \ \textnormal{ and } \ \ \Hom[\Ropp]{\bi{I}}{\bi{I}}.
\end{displaymath}
Moreover, the $S$--linearity of $\lambda$ makes the above diagram
commutative. To see this, observe that $(\Hom[\Ropp]{\lambda}{\bi{I}}
\circ \lhty[I])(s)$ and $(\Hom[\Ropp]{\bi{P}}{\lambda} \circ
\lhty[P])(s)$ yield maps
\begin{gather*}
  \xymatrix{ {}_S P_R \ar[r]^-{\lambda} & {}_S I_R \ar[r]^-{s \cdot -}
    & {}_S I_R } \quad\text{and}\quad \xymatrix{ {}_S P_R \ar[r]^-{s
      \cdot -} & {}_S P_R \ar[r]^-{\lambda} & {}_S I_R , }
\end{gather*}
respectively, where \mbox{$s\cdot-$} denotes left-multiplication with a generic
element in $S$. A similar analysis shows that the $S$--linearity of
$\lambda$ implies both $S$--linearity of
$\Hom[\Ropp]{\bi{P}}{\lambda}$ and $\Sopp$--linearity of
$\Hom[\Ropp]{\lambda}{\bi{I}}$.

Since $\bi{P}$ is a right-bounded complex of projective
$\Ropp$--modules, $\Hom[\Ropp]{\bi{P}}{\lambda}$ is a
quasi-isomorphism.  Similarly, $\Hom[\Ropp]{\lambda}{\bi{I}}$ is a
quasi-isomorphism.  Consequently,
\begin{displaymath}
  \lhty[P] \textnormal{ is a quasi-isomorphism } 
  \iff 
  \lhty[I] \textnormal{ is a quasi-isomorphism}.
\end{displaymath}
When we in \eqref{dfn:dc}(4) require that $\lhty \colon _SS_S
\longrightarrow \DHom[\Ropp]{\bi{D}}{\bi{D}}$ is invertible in $\D[S]$
(equivalently, invertible in $\D[\Sopp]$), it means that $\lhty[P]$ is
a quasi-isomorphism of $S$--complexes (equivalently, of
$\Sopp$--complexes).

Similar remarks apply to the morphism $\rhty \colon _RR_R
\longrightarrow \DHom[S]{\bi{D}}{\bi{D}}$; here it becomes important
that $\lambda$ is $\Ropp$--linear.

\begin{proof}[Proof of \eqref{prp:OppPair}]
  \label{proof:OppPair}
  By symmetry it suffices to prove that a dualizing complex $D$ for a
  noetherian pair $\pair$ is dualizing for $\pair[\Ropp,\Sopp]$ as
  well.  Obviously, $\Ropp$ is left noetherian and $\Sopp$ is right
  noetherian, so $\pair[\Ropp,\Sopp]$ is a noetherian pair.
  Furthermore, since $(S,\Ropp)$--bimodules are naturally identified
  with $(\Ropp,(\Sopp)^{\mathrm{opp}})$--bimodules, it is clear that
  $D$ satisfies conditions (1), (2), and (3) in \eqref{dfn:dc} relative
  to $\pair[\Ropp,\Sopp]$. Finally, we need to see that the homothety
  morphisms are invertible. The morphism
  \begin{displaymath}
    \acute{\chi}_P^{\pair[\Ropp,\Sopp]} \colon \Sopp \longrightarrow
    \Hom[\Ropp]{P}{P} 
  \end{displaymath}
  is identical to $\mapdef{\lhty[P]}{S}{\Hom[\Ropp]{P}{P}}$ through
  the identification $S=\Sopp$ (as $(S,\Sopp)$--bimodules, not as
  rings). By assumption $\lhty[P]$ is a quasi-isomorphism, and hence
  so is $\acute{\chi}_P^{\pair[\Ropp,\Sopp]}$. Similarly,
  $\grave{\chi}_P^{\pair[\Ropp,\Sopp]}$ is identified with
  $\rhty[P]$.
\end{proof}

\begin{ipg}
  The next result was also stated in section 1 (prop.\ 
  \eqref{prp:PJresult}); this time we prove it.
\end{ipg}

\begin{prp}
  \label{bfhpg:PJproof}
  Assume that the noetherian pair $\pair$ has a dualizing complex
  $\bi{D}$. If $X \in \D$ has finite $\fd{X}$, then there is an
  inequality,
  \begin{displaymath}  
    \pd{X} \le \max\big\{\id[S]{\biP{D}} + \supP{\Dtp{\bi{D}}{X}} \,,\,
    \sup{X}\big\} < \infty. 
  \end{displaymath}
  Moreover, $\FPD$ is finite if and only if $\,\FFD$ is finite.
\end{prp} 

\begin{proof}
  Define the integer $n$ by
  \begin{equation*}
    n = \max\big\{\id[S]{\biP{D}} + \supP{\Dtp{\bi{D}}{X}} \,,\,
    \sup{X}\big\} 
  \end{equation*}
  Let $Q \xre X$ be a projective resolution of $X$. Since $\infty > n
  \ge \sup{X} = \sup{Q}$, we have a quasi-isomorphism $Q \xre
  \Tsl{n}{Q}$, and hence it suffices to prove that the $R$--module
  $\cC{n}{Q} = \Coker{(Q_{n+1} \lora Q_n)}$ is projective. This is
  tantamount to showing that $\Ext{1}{\cC{n}{Q}}{\cC{n+1}{Q}}=0$.  We
  have the following isomorphisms of abelian groups:
  \begin{align*}
    \Ext{1}{\cC{n}{Q}}{\cC{n+1}{Q}} &\,\cong\,
    \operatorname{H}_{-(n+1)}
    \DHom{X}{\cC{n+1}{Q}} \\
    &\,\cong\, \operatorname{H}_{-(n+1)}
    \DHom{X}{\DHom[S]{\bi{D}}{\Dtp{\bi{D}}{\cC{n+1}{Q}}}} \\
    &\,\cong\, \operatorname{H}_{-(n+1)}
    \DHom[S]{\Dtp{\bi{D}}{X}}{\Dtp{\bi{D}}{\cC{n+1}{Q}}}
  \end{align*}
  The first isomorphism follows as $n \ge \sup{X}$, and the second one
  follows as $\fd{\cC{n+1}{Q}}$ is finite, and hence $\cC{n+1}{Q} \in
  \A$ by \eqref{bfhpg:Foxby_eq}. The third isomorphism is by
  adjointness.  It is now sufficient to show that
  \begin{displaymath}
    -\inf \DHom[S]{\Dtp{\bi{D}}{X}}{\Dtp{\bi{D}}{\cC{n+1}{Q}}}
    \leq n,
  \end{displaymath}
  and this follows as:
  \begin{gather*}
    -\inf \DHom[S]{\Dtp{\bi{D}}{X}}{\Dtp{\bi{D}}{\cC{n+1}{Q}}}\qquad\qquad\qquad\qquad\\
    \qquad\qquad\qquad\qquad
    \begin{split}
      &\le\, \id[S]{\DtpP{\bi{D}}{\cC{n+1}{Q}}} + \supP{\Dtp{\bi{D}}{X}}\\
      &\le\, \id[S]{\biP{D}} + \supP{\Dtp{\bi{D}}{X}}\\
      &\le\, n.
    \end{split}
  \end{gather*}
  The first inequality is by \cite[thm.~2.4.I]{LLAHBF91}, and the
  second one is by \cite[thm.~4.5(F)]{LLAHBF91}, as $S$ is left
  noetherian and $\fd{\cC{n+1}{Q}}$ is finite.
  
  The last claim follows as
  \begin{displaymath}
    \FFD \,\le\, \FPD \,\le\, \FFD + \id[S]{\biP{D}} + \supP{\bi{D}}.
  \end{displaymath}
  The first inequality is by \cite[prop.~6]{CUJ70}. To verify the
  second one, let $M$ be a module with $\pd{M}$ finite. We have
  already seen that
  \begin{displaymath}
    \pd{M} \le \max\big\{\id[S]{\biP{D}} + \supP{\Dtp{\bi{D}}{M}} \,,\,
    0\big\},
  \end{displaymath}
  so it suffices to see that
  \begin{equation*}
    \supP{\Dtp{\bi{D}}{M}} \le \fd{M} + \supP{\bi{D}},
  \end{equation*}
  and that follows from \cite[thm.~2.4.F]{LLAHBF91}.
\end{proof}

%%% SECTION A

\section*{Acknowledgments}
It is a pleasure to thank Srikanth Iyengar, Peter J{\o}rgensen, Sean
Sather-Wagstaff and Oana Veliche for their readiness to discuss the
present work as well as their own work in this field. This paper has
benefited greatly from their suggestions. Finally, we are very
grateful for the referee's good and thorough comments.

\bibliographystyle{amsplain} \providecommand{\bysame}{\leavevmode\hbox

\begin{thebibliography}{10}

  
\bibitem{ALL} Leovigildo Alonso~Tarr{\'{\i}}o, Ana
  Jerem{\'{\i}}as~L{\'o}pez, and Joseph Lipman, \emph{Local homology
    and cohomology on schemes}, Ann. Sci. \'Ecole Norm. Sup. (4)
  \textbf{30} (1997), no.~1, 1--39. \MR{98d:14028}
  
\bibitem{MAu67} Maurice Auslander, \emph{Anneaux de {G}orenstein, et
    torsion en alg\`ebre commutative}, Secr\'etariat math\'ematique,
  Paris, 1967, S\'eminaire d'Alg\`ebre Commutative dirig\'e par Pierre
  Samuel, 1966/67. Texte r\'edig\'e, d'apr\`es des expos\'es de
  Maurice Auslander, par Marquerite Mangeney, Christian Peskine et
  Lucien Szpiro. \'Ecole Normale Sup\'erieure de Jeunes Filles.
  
\bibitem{MAuMBr69} Maurice Auslander and Mark Bridger, \emph{Stable
    module theory}, Memoirs of the American Mathematical Society, No.
  94, American Mathematical Society, Providence, R.I., 1969. \MR{42
    \#4580}
  
\bibitem{LLAHBF91} Luchezar~L. Avramov and Hans-Bj{\o}rn Foxby,
  \emph{Homological dimensions of unbounded complexes}, J. Pure Appl.
  Algebra \textbf{71} (1991), no.~2-3, 129--155. \MR{93g:18017}
  
\bibitem{LLAHBF97} \bysame, \emph{Ring homomorphisms and finite
    {G}orenstein dimension}, Proc.  London Math. Soc. (3) \textbf{75}
  (1997), no.~2, 241--270. \MR{98d:13014}
  
\bibitem{LLAHBFBHr94} Luchezar~L. Avramov, Hans-Bj{\o}rn Foxby, and
  Bernd Herzog, \emph{Structure of local homomorphisms}, J. Algebra
  \textbf{164} (1994), no.~1, 124--145.  \MR{95f:13029}
  
\bibitem{HBs62} Hyman Bass, \emph{Injective dimension in {N}oetherian
    rings}, Trans. Amer.  Math. Soc. \textbf{102} (1962), 18--29.
  \MR{25 \#2087}
  
\bibitem{LCh76} Leo~G. Chouinard, II, \emph{On finite weak and
    injective dimension}, Proc.  Amer. Math. Soc. \textbf{60} (1976),
  57--60 (1977). \MR{54 \#5217}
  
\bibitem{LWC3} Lars~Winther Christensen, \emph{Gorenstein dimensions},
  Lecture Notes in Mathematics, vol. 1747, Springer-Verlag, Berlin,
  2000. \MR{2002e:13032}
  
\bibitem{LWC1} \bysame, \emph{Semi-dualizing complexes and their
    {A}uslander categories}, Trans. Amer. Math. Soc. \textbf{353}
  (2001), no.~5, 1839--1883 (electronic).  \MR{2002a:13017}
  
\bibitem{RHD} Lars~Winther Christensen, Hans-Bj{\o}rn Foxby, and
  Anders Frankild, \emph{Restricted homological dimensions and
    {C}ohen-{M}acaulayness}, J.  Algebra \textbf{251} (2002), no.~1,
  479--502. \MR{2003e:13022}
  
\bibitem{LWCHH} Lars~Winther Christensen and Henrik Holm, \emph{Ascent
    properties of {A}uslander categories}, preprint (2005), available
  from \mbox{\sffamily http://arXiv.org/math.AC/0509570}.
  
\bibitem{EEn81} Edgar~E. Enochs, \emph{Injective and flat covers,
    envelopes and resolvents}, Israel J. Math. \textbf{39} (1981),
  no.~3, 189--209. \MR{83a:16031}
  
\bibitem{EEnOJn93} Edgar~E. Enochs and Overtoun M.~G. Jenda, \emph{On
    {G}orenstein injective modules}, Comm. Algebra \textbf{21} (1993),
  no.~10, 3489--3501.  \MR{94g:13006}
  
\bibitem{EEnOJn95} \bysame, \emph{Gorenstein injective and projective
    modules}, Math. Z.  \textbf{220} (1995), no.~4, 611--633.
  \MR{97c:16011}
  
\bibitem{EEnOJn95a} \bysame, \emph{Resolutions by {G}orenstein
    injective and projective modules and modules of finite injective
    dimension over {G}orenstein rings}, Comm. Algebra \textbf{23}
  (1995), no.~3, 869--877. \MR{96m:16010}
  
\bibitem{EJT93} Edgar~E. Enochs, Overtoun M.~G. Jenda, and Blas
  Torrecillas, \emph{Gorenstein flat modules}, Nanjing Daxue Xuebao
  Shuxue Bannian Kan \textbf{10} (1993), no.~1, 1--9. \MR{95a:16004}
  
\bibitem{EJX96} Edgar~E. Enochs, Overtoun M.~G. Jenda, and Jin~Zhong
  Xu, \emph{Foxby duality and {G}orenstein injective and projective
    modules}, Trans. Amer. Math. Soc.  \textbf{348} (1996), no.~8,
  3223--3234. \MR{96k:13010}
  
\bibitem{EEnJAL02} Edgar~E. Enochs and J.~A. L{\'o}pez-Ramos,
  \emph{Kaplansky classes}, Rend. Sem.  Mat. Univ. Padova \textbf{107}
  (2002), 67--79. \MR{MR1926201 (2003j:16005)}
  
\bibitem{HHA} Hans-Bj{\o}rn Foxby, \emph{Hyperhomological algebra \&
    commutative rings}, notes in preparation.
  
\bibitem{HBF77} \bysame, \emph{Isomorphisms between complexes with
    applications to the homological theory of modules}, Math. Scand.
  \textbf{40} (1977), no.~1, 5--19. \MR{56 \#5584}
  
\bibitem{HBF79} \bysame, \emph{Bounded complexes of flat modules}, J.
  Pure Appl. Algebra \textbf{15} (1979), no.~2, 149--172.
  \MR{83c:13008}
  
\bibitem{HBF94} \bysame, \emph{Gorenstein dimension over
    {C}ohen-{M}acaulay rings}, Proceedings of international conference
  on commutative algebra (W.~Bruns, ed.), Universit\"at Onsabr\"uck,
  1994.
  
\bibitem{SriHBF} Hans-Bj{\o}rn Foxby and Srikanth Iyengar, \emph{Depth
    and amplitude for unbounded complexes}, Commutative algebra.
  Interaction with Algebraic Geometry (Grenoble-Lyon 2001), Amer.
  Math. Soc., Providence, RI, no. 331, 2003, pp.~119--137.
  
\bibitem{F} Anders Frankild, \emph{Vanishing of local homology}, Math.
  Z. \textbf{244} (2003), no.~3, 615--630. \MR{1 992 028}
  
\bibitem{AF-SI-PJ:DDGMGDGA} Anders Frankild, Srikanth Iyengar, and
  Peter J{\o}rgensen, \emph{Dualizing differential graded modules and
    {G}orenstein differential graded algebras}, J. London Math. Soc.
  (2) \textbf{68} (2003), no.~2, 288--306. \MR{MR1994683
    (2004f:16013)}
  
\bibitem{GM} Sergei~I. Gelfand and Yuri~I. Manin, \emph{Methods of
    homological algebra}, second ed., Springer Monographs in
  Mathematics, Springer-Verlag, Berlin, 2003. \MR{2003m:18001}
  
\bibitem{GreenMay} J.~P.~C. Greenlees and J.~P. May, \emph{Derived
    functors of {$I$}-adic completion and local homology}, J. Algebra
  \textbf{149} (1992), no.~2, 438--453. \MR{93h:13009}
  
\bibitem{egaIII} A.~Grothendieck, \emph{\'{E}l\'ements de
    g\'eom\'etrie alg\'ebrique. {III}.  \'{E}tude cohomologique des
    faisceaux coh\'erents. {I}}, Inst. Hautes \'Etudes Sci. Publ.
  Math. (1961), no.~11, 167. \MR{MR0163910 (29 \#1209)}
  
\bibitem{RAD} Robin Hartshorne, \emph{Residues and duality}, Lecture
  notes of a seminar on the work of A. Grothendieck, given at Harvard
  1963/64. With an appendix by P.  Deligne. Lecture Notes in
  Mathematics, No. 20, Springer-Verlag, Berlin, 1966.  \MR{36 \#5145}
  
\bibitem{HH3} Henrik Holm, \emph{Gorenstein homological dimensions},
  J. Pure Appl. Algebra \textbf{189} (2004), no.~1-3, 167--193.
  \MR{2038564 (2004k:16013)}
  
\bibitem{Ishikawa} Ishikawa, Takeshi, \emph{On injective modules and
    flat modules}, J. Math. Soc. Japan \textbf{17} (1965), 291--296.
  \MR{0188272 (32 \#5711)}
  
\bibitem{SriSean} Srikanth Iyengar and Sean Sather-Wagstaff,
  \emph{G-dimension over local homomorphisms. {A}pplications to the
    {F}robenius endomorphism}, Illinois J.  Math. \textbf{48} (2004),
  no.~1, 241--272. \MR{2048224}
  
\bibitem{CUJ70} C.~U. Jensen, \emph{On the vanishing of
    {$\underset{\longleftarrow}{\lim}^{(i)}$}}, J. Algebra \textbf{15}
  (1970), 151--166. \MR{41 \#5460}
  
\bibitem{DAJLMS} David~A. Jorgensen and Liana~M. \c{S}ega,
  \emph{Independence of the total reflexivity conditions for modules},
  to appear in {A}lgebr.\ {R}epresent.\ {T}heory, preprint (2004)
  available from \mbox{\sffamily http://arXiv.org/math.AC/0410257}.
  
\bibitem{PJ:Preprint} Peter J{\o}rgensen, \emph{Finite flat and
    projective dimension}, Comm. Algebra \textbf{33} (2005), no.~7,
  2275--2279. \MR{MR2153221}
  
\bibitem{TKw02} Takesi Kawasaki, \emph{On arithmetic
    {M}acaulayfication of {N}oetherian rings}, Trans. Amer. Math. Soc.
  \textbf{354} (2002), no.~1, 123--149 (electronic).  \MR{2002i:13001}
  
\bibitem{DLz69} Daniel Lazard, \emph{Autour de la platitude}, Bull.
  Soc. Math. France \textbf{97} (1969), 81--128. \MR{40:7310}
  
\bibitem{EMt74} Eben Matlis, \emph{The {K}oszul complex and duality},
  Comm. Algebra \textbf{1} (1974), 87--144. \MR{MR0344241 (49 \#8980)}
  
\bibitem{Mat} Hideyuki Matsumura, \emph{Commutative ring theory},
  second ed., Cambridge Studies in Advanced Mathematics, vol.~8,
  Cambridge University Press, Cambridge, 1989, Translated from the
  Japanese by M. Reid. \MR{90i:13001}
  
\bibitem{Miyachi:DCMDT} Jun-ichi Miyachi, \emph{Derived categories and
    {M}orita duality theory}, J.  Pure Appl. Algebra \textbf{128}
  (1998), no.~2, 153--170. \MR{MR1624752 (99d:16004)}
  
\bibitem{MRnLGr71} Michel Raynaud and Laurent Gruson, \emph{Crit\`eres
    de platitude et de projectivit\'e. {T}echniques de
    ``platification'' d'un module}, Invent. Math.  \textbf{13} (1971),
  1--89. \MR{46 \#7219}
  
\bibitem{PSc03} Peter Schenzel, \emph{Proregular sequences, local
    cohomology, and completion}, Math. Scand. \textbf{92} (2003),
  no.~2, 161--180. \MR{MR1973941 (2004f:13023)}
  
\bibitem{MLT76} Mark~L. Teply, \emph{Torsion-free covers. {II}},
  Israel J. Math. \textbf{23} (1976), no.~2, 132--136. \MR{MR0417245
    (54 \#5302)}
  
\bibitem{OV03} Oana Veliche, \emph{Gorenstein projective dimension for
    complexes}, Trans.  Amer. Math. Soc. \textbf{358} (2006),
  1257--1283.

\bibitem{Verdier} Jean-Louis Verdier, \emph{Des cat\'egories
    d\'eriv\'ees des cat\'egories ab\'eliennes}, Ast\'erisque (1996),
  no.~239, xii+253 pp. (1997), With a preface by Luc Illusie, Edited
  and with a note by Georges Maltsiniotis.  \MR{98c:18007}
  
\bibitem{Wei} Charles A. Weibel, \emph{An introduction to homological
    algebra}, Cambridge studies in advanced mathematics, vol.~38,
  Cambridge University Press, Cambridge, 1994, xiv+450.
  
\bibitem{Wu-Zhang:DCncLR} Q.-S. Wu and J.~J. Zhang, \emph{Dualizing
    complexes over noncommutative local rings}, J. Algebra
  \textbf{239} (2001), no.~2, 513--548. \MR{MR1832904 (2002g:16010)}
  
\bibitem{SYs95} Siamak Yassemi, \emph{G-dimension}, Math. Scand.
  \textbf{77} (1995), no.~2, 161--174.
  
\bibitem{SYs98} \bysame, \emph{Width of complexes of modules}, Acta
  Math. Vietnam. \textbf{23} (1998), no.~1, 161--169. \MR{99g:13026}
  
\bibitem{Yekutieli:DCncGA} Amnon Yekutieli, \emph{Dualizing complexes
    over noncommutative graded algebras}, J. Algebra \textbf{153}
  (1992), no.~1, 41--84. \MR{MR1195406 (94a:16077)}
  
\bibitem{Yekutieli-Zhang:RADC} Amnon Yekutieli and James~J. Zhang,
  \emph{Rings with {A}uslander dualizing complexes}, J. Algebra
  \textbf{213} (1999), no.~1, 1--51. \MR{MR1674648 (2000f:16012)}

\end{thebibliography}
  to3em{\hrulefill}\thinspace}
\providecommand{\MR}{\relax\ifhmode\unskip\space\fi MR }

\providecommand{\MRhref}[2]{%
  \href{http://www.ams.org/mathscinet-getitem?mr=#1}{#2} }
\providecommand{\href}[2]{#2}

\end{document}